\documentclass[11pt,leqno]{amsart}
\usepackage{amsthm,amsfonts,amssymb,amsmath,oldgerm,mathrsfs,mathabx}
\usepackage[scr=boondoxo]{mathalfa}
\numberwithin{equation}{section}
\usepackage{fullpage,setspace,fancyhdr,bm}
\usepackage{graphicx,psfrag,listings} 
\usepackage{color}

\usepackage[top=30mm,bottom=30mm,left=25mm,right=25mm,a4paper]{geometry}

\usepackage{hyperref}
\usepackage{graphics}

\newtheorem{definition}{Definition}
\newtheorem{lemma}{Lemma}
\newtheorem{corollary}{Corollary}
\newtheorem{proposition}{Proposition}
\newtheorem{theorem}{Theorem}
\newtheorem{remark}{Remark}

\newtheorem{assumption}{Assumption}

\renewcommand\d{\partial}

\DeclareMathOperator{\iD}{i}
\DeclareMathOperator{\eD}{e}
\DeclareMathOperator{\dD}{d}

\newcommand{\Id}{\text{Id}}
\newcommand{\I}{\text{I}}

\newcommand{\One}{\text{\bf 1}}

\DeclareMathOperator{\Div}{div}

\newcommand{\fpar}{f_{\shortparallel}}
\newcommand{\fperp}{\bF_{\perp}}
\newcommand{\sigpar}{\sigma_{\shortparallel}}
\newcommand{\sigperp}{\bfsigma_{\perp}}

\newcommand\br{\begin{remark}}
\newcommand\er{\end{remark}}
\newcommand\bp{\begin{pmatrix}}
\newcommand\ep{\end{pmatrix}}
\newcommand{\be}{\begin{equation}}
\newcommand{\ee}{\end{equation}}

\newcommand{\bpr}{\begin{proposition}}
\newcommand{\epr}{\end{proposition}}
\newcommand{\bt}{\begin{theorem}}
\newcommand{\et}{\end{theorem}}
\newcommand{\bc}{\begin{corollary}}
\newcommand{\ec}{\end{corollary}}
\newcommand{\bl}{\begin{lemma}}
\newcommand{\el}{\end{lemma}}
\newcommand{\bd}{\begin{definition}}
\newcommand{\ed}{\end{definition}}

\newcommand\R{\mathbf R}
\newcommand\C{\mathbf C}
\newcommand{\N}{\mathbf N}

\newcommand\bF{{\mathbf F}}

\newcommand\bS{{\mathbf S}}

\newcommand\bY{{\mathbf Y}}
\newcommand\bZ{{\mathbf Z}}

\newcommand\bfx{{\mathbf x}}
\newcommand\bfy{{\mathbf y}}

\newcommand\bfeta{{\boldsymbol \eta}}

\newcommand\bfsigma{{\boldsymbol \sigma}}

\newcommand\uu{{\underline u}}

\newcommand\cC{{\mathcal C}}
\newcommand\cD{{\mathcal D}}

\newcommand\cT{{\mathcal T}}
\newcommand\cU{{\mathcal U}}

\newcommand\cY{{\mathcal Y}}

\newcommand*{\front}{\uu}

\title{Convective stability in scalar balance laws}

\author{Louis Gar\'enaux}
\address{Karlsruhe Institute for Technology, Englerstraße 2, 76131 Karlsruhe,
Germany}
\email{{\tt louis.garenaux@kit.edu}}
\thanks{}

\author{L.~Miguel Rodrigues}
\address{
Univ Rennes, CNRS, IRMAR - UMR 6625, F-35000 Rennes, France}
\email{{\tt luis-miguel.rodrigues@univ-rennes.fr}}
\thanks{Research of L.M.R. was supported by the Institut Universitaire de France and partially supported by France 2030 through the programme Centre Henri Lebesgue ANR-11-LABX-0020-01.
Research of L.G. was funded by the Deutsche Forschungsgemeinschaft (DFG, German Research Foundation) -- Project-ID 258734477 -- SFB 1173.
}

\begin{document}

\maketitle

\begin{abstract}
We fill the two main remaining gaps in the full classification of non-degenerate planar traveling waves of scalar balance laws from the point of view of spectral and nonlinear stability/instability under smooth perturbations. On one hand we investigate the impact on the classification of allowing to restrict perturbations by further localization constraints, the so-called \emph{convective} stability analysis. On the other hand, in the multidimensional case we extend the existing analysis to general perturbations, relaxing the artificial restriction, assumed so far, that initial perturbations are supported away from characteristic points. A striking feature of the latter analysis is that it involves infinite-dimensional families of genuinely multidimensional traveling waves that live near the perturbed planar traveling wave.

\vspace{0.5em}

{\small \paragraph {\bf Keywords:} discontinuous traveling waves; characteristic points; convective stability; scalar balance laws.
}

\vspace{0.5em}

{\small \paragraph {\bf AMS Subject Classifications:} 35B35, 35L02, 35L67, 35B40, 35L03, 35P15, 37L15.
}
\end{abstract}

\tableofcontents

\section{Introduction}

In the present contribution, we continue the study, initiated in~\cite{DR1, DR2}, of the large-time asymptotic behavior of solutions to first order scalar hyperbolic balance laws, in neighborhoods of traveling waves. 

\subsection{Brief overview of results}

Explicitly, let $f$ and $g$ be two smooth scalar functions on $\R$ and consider the associated one-dimensional scalar balance law
\begin{equation}
\label{e:original-eq}
\partial_t u + \partial_x (f(u)) = g(u)
\end{equation}
for the unknown function $u$, $u(t,x)\in\R$ depending on time variable $t$ and spatial variable $x\in\R$. A significant part of our contribution is also devoted to planar traveling waves of \emph{multidimensional} balance laws but, for the sake of exposition, we mostly restrict the introduction to the one-dimensional equation \eqref{e:original-eq} and delay to Section~\ref{s:multiD} the multidimentional discussion.

A \emph{traveling-wave} solution to \eqref{e:original-eq} is a solution of the form $u(t,x) = \front(x - \sigma t)$, for some speed $\sigma$ and some profile $\front$. In other words traveling waves are equilibria in some frame moving at a constant speed. Note that our definition includes as special cases, standing waves that correspond to the case when $\sigma=0$, and constant equilibria that correspond to the case when the profile $\front$ is constant. Our goal is to understand the dynamics near such special solutions.  Obviously, as relative equilibria, traveling waves are expected to play a prominent role in the description of the large-time dynamics. For general background on traveling waves and their stability we refer to \cite{KapitulaPromislow-stability}.

As in~\cite{DR1, DR2}, we allow profiles to be discontinuous, but, then, restrict, among possible weak solutions, to entropy-admissible solutions. In particular, the classical Kru\v zkov theory~\cite{Kruzhkov_70} applies to the problems we study, and yields uniqueness for the notion of solution we use throughout. The classical motivation for the consideration of non smooth solutions is that even smooth initial data may lead to the formation of discontinuities in finite time. For more on the latter and other relevant elementary background on hyperbolic equations, we refer the reader to~\cite{Bressan_00}.

A dramatic consequence of allowing discontinuous solutions is that Equation~\eqref{e:original-eq} may possess a tremendously huge number of traveling-wave solutions. Yet one of the main upshots of the analysis of \cite{DR2} is that very few of those are stable. Let us describe the analysis underlying the claim in a few words. In the framework of piecewise-smooth \emph{unweighted} topologies \cite{DR2} provides
\begin{itemize}
\item a complete classification of non-degenerate traveling waves according to their spectral stability;
\item proofs that for those waves spectral instability (resp. stability) yields dynamical nonlinear instability (resp. asymptotic stability).
\end{itemize}
The classification relies on the fact, proved in \cite{DR2}, that, in such a framework, spectral (and nonlinear) instability occur if and only if the profile exhibits at least one of the following features:
\begin{itemize}
\item an endstate $u_\infty$ --- that is, a limit of $\front$ at $+\infty$ or $-\infty$ --- such that $g'(u_\infty)>0$;
\item a discontinuity point $d_0$ at which $\frac{[\,g(\front)\,]_{d_0}}{[\,\front\,]_{d_0}}>0$;
\item a characteristic value $u_\star$ --- that is, a value $u_\star$ of $\front$ with $f'(u_\star)=\sigma$ --- such that $g'(u_\star)<0$. 
\end{itemize}
In the above, we have used notation $[A]_{d}$ to denote the jump of a function $A$ at a point $d$,
\begin{align*}
[A]_d&:=A(d^+)-A(d^-)\,,&
A(d^+)&:=\lim_{\substack{x\to d\\x>d}}A(x)\,,&
A(d^-)&:=\lim_{\substack{x\to d\\x<d}}A(x)\,.
\end{align*}

Extending the foregoing analysis, our present one-dimensional achievement is the elucidation of which of the instabilities due to an unstable endstate may be stabilized by working in a \emph{weighted} topology. In particular we describe the dynamics near monostable or even nullistable fronts -- i.e. structures with respectively exactly $1$ or $0$ asymptotically stable endstate. A general motivation for such kind of studies is that even though those waves are unstable in unweighted topologies they may play the role of an elementary block in the description of the large-time dynamics of compactly supported\footnote{In the sense that on each connected component of a neighborhood of infinity they coincide with a constant solution of the equation.} initial data. See for instance \cite{Sattinger_76, Aronson_Weinberger_78} for the most classical studies in a reaction-diffusion setting. We expect the case of Equation~\eqref{e:original-eq} to be no exception. In this direction, we mention that, under rather stringent but natural assumptions on $f$ and $g$ --- including the strict convexity of $f$ and the strict dissipativity at infinity of $g$ ---, it has already been proved in \cite{Sinestrari_96,Mascia_Sinestrari_97}, by comparison principle arguments, that starting from an $L^\infty$ initial datum constant near $-\infty$ and near $+\infty$, the large-time dynamics is well captured in $L^\infty$ topology by piecing together traveling waves (constants, fronts or periodic waves). We also refer the reader to \cite{RYZ,FRYZ} for some supporting numerical experiments and, accompanying spectral studies, in some system cases.

At the linearized level, all the instabilities due to a non-characteristic\footnote{That is, for the limiting endstate $u_\infty$ of a wave at speed $\sigma$, such that $f'(u_\infty)\neq \sigma$.} unstable endstate may be stabilized with a suitably designed weight. Yet only the cases when the spectral stabilization may be carried out with a weight enforcing decay at infinity are directly relevant for the nonlinear stability analysis. One conclusion of our analysis is that, as may be guessed from heuristic arguments, endstate instabilities may be nonlinearly stabilized provided that unstable wave-packets, or, equivalently here, characteristics, travel outward the unstable infinity. Note in particular that unstable constant states cannot be stabilized in this way, since they are unstable at both infinities but with characteristics pointing outward at one infinity, and inward at the other one. It is precisely because of this connection between the motion of unstable wave-packets and stability in weighted spaces that the latter is often referred to as \emph{convective} stability. 

With this in mind, we stress that there is little hope to stabilize the other kinds of instability identified in \cite{DR2} since in remaining cases the unstable wave-packets point towards the instability sources (either characteristic points or points of discontinuity). It goes without writing that all the corresponding instability results from \cite{DR2} hold essentially unchanged in the kind of weighted spaces considered here, that only modify localization at infinities.

\begin{figure}
\begin{center}
 \begin{tabular}{cc}
\includegraphics[width=7.75cm]{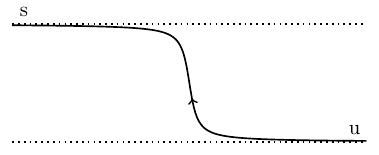} &    
\includegraphics[width=7.75cm]{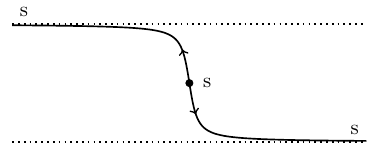} 
\\
(a)  & (b)  \\
\includegraphics[width=7.75cm]{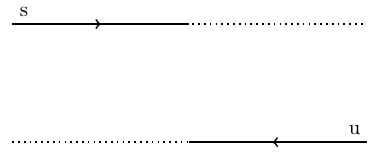} &    
\includegraphics[width=7.75cm]{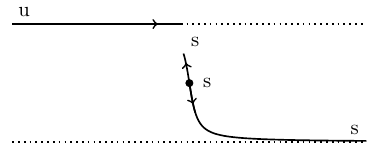} 
\\
(c)  & (d)  
\end{tabular}
\caption{\label{fig:profiles}
{\bf Wave profiles.} A few shapes for wave profiles with (a) a smooth non-characteristic monostable front (b) a smooth bistable characteristic front (c) a monostable Riemann shock and (d) a characteristic discontinuous monostable front. Arrows show the direction of propagation of characteristics, dots stand for characteristic points, letters $u$ and $s$ indicate whether the endstate is stable or unstable. Entropy conditions determine the direction of characteristics near discontinuities. For non constant parts stability/instability of the endstate determine the direction of characteristics near infinity.}
 \end{center}
\end{figure}

In contrast, it follows from the wave profile equation itself that if $u_\infty$ is a non-characteristic endstate at $-\infty$ (respectively $+\infty$) that is non degenerate in the sense that $g'(u_\infty)\neq0$ and the wave profile is not constant near the infinity at hand then $g'(u_\infty)/(f'(u_\infty)-\sigma)$ is positive (resp. negative). As a consequence, when $u_\infty$ is unstable, that is, $g'(u_\infty)>0$, and the corresponding portion of the wave profile is not constant, characteristics are outgoing from the infinity at hand. Besides, on portions of wave profiles that are constant between a discontinuity and one infinity, independently of any stability consideration the entropy conditions yield the same conclusion, characteristics move outward infinity, toward the discontinuity. Consistently we do prove that when dealing with non constant\footnote{We allow however constant portions.} traveling waves all degenerate instabilities due to endstate instabilities may be stabilized by using spatial weights. We illustrate this discussion in Figure~\ref{fig:profiles}.

A particular interesting consequence is that though non degenerate smooth non characteristic waves are all unstable in unweighted topologies, they are all convectively stable. Let us give a more concrete description of such waves. Pick two consecutive zeros of $g$, $u_{-\infty}$ and $u_{+\infty}$ (without assuming any particular order\footnote{In particular, in this paragraph by convention $[u_{-\infty},u_{+\infty}]=[\min(\{u_{-\infty},u_{+\infty}\}),\max(\{u_{-\infty},u_{+\infty}\})]$.}). Assume that they are non degenerate as zeros of $g$, so that $g'(u_{-\infty})\,g'(u_{+\infty})<0$. Then if $g'(u_{-\infty})>0$ (respectively $g'(u_{+\infty})>0$), for any $\sigma<\min_{[u_{-\infty},u_{+\infty}]}(f')$ (resp. $\sigma>\max_{[u_{-\infty},u_{+\infty}]}(f')$), there exists a smooth traveling wave at speed $\sigma$ connecting $u_{-\infty}$ at $-\infty$ to $u_{+\infty}$ at $+\infty$, unique up to spatial translations. Note that the sign conditions on $f'-\sigma$ are consistent with the intuition that the stable endstate invades the unstable one. This type of monostable front has been the object of considerable studies in scalar reaction-diffusion equations. Our analysis of the present advection-reaction equation shares many common points with those. Yet a crucial point of divergence is that in the present case there is no traveling wave of critical speed ($\min_{[u_{-\infty},u_{+\infty}]}(f')$ if $g'(u_{-\infty})>0$, $\max_{[u_{-\infty},u_{+\infty}]}(f')$ if $g'(u_{+\infty})>0$). This is all the more striking that in reaction-diffusion equations critical speed waves are precisely the elementary blocks involved in the description of the large-time dynamics arising from compactly supported initial data \cite{Bramson_83,Roquejoffre-RoussierMichon}.

\subsection{A few technical insights}

To provide an introduction to the technical side of our analysis, we first make more precise some terminology. 
\bd\label{def:non-degenerate}
We say that a wave $(t,x)\mapsto \uu(x-\sigma\,t)$ is non-degenerate if its wave profile $\uu$ is piecewise smooth (with a discrete set of discontinuity points), bounded, and such that 
\begin{itemize}
\item any characteristic value of the profile $u_\star$ satisfies $f''(u_\star)\neq 0$;
\item at any discontinuity point $d$, none of the limiting values $\uu(d^-)$ and $\uu(d^+)$ is a characteristic value and Oleinik's conditions are satisfied in a strict sense:\\ $f'(\uu(d^-))>\sigma>f'(\uu(d^+))$ and for any $v$ between $\uu(d^-)$ and $\uu(d^+)$,
\[
\frac{f'(v)-f'(\uu(d^-))}{v-\uu(d^-)}>\frac{f'(v)-f'(\uu(d^+))}{v-\uu(d^+)}\,;
\]
\item any endstate $u_\infty$ satisfies $g'(u_\infty)\neq0$ and is not a characteristic value.
\end{itemize}
\ed
The Oleinik conditions are equivalent to entropy admissibility for the wave solution. Enforcing them in strict form ensures that nearby weak solutions are also entropy admissible. Thus, as already pointed out, a form of well-posedness is obtained from the seminal results of Kru\v zkov \cite{Kruzhkov_70}. 

Our strategy of proof for such waves follows closely the one designed in \cite{DR1} and subsequently used in \cite{DR2}. In particular, we treat discontinuities indirectly. First we show the stability of smooth blocks, then we extend each smooth portion so as to be able to apply the latter, and finally  we glue the smooth blocks along discontinuity positions determined from the Rankine-Hugoniot conditions. The latter conditions encode that by gluing smooth blocks solving the equation in a classical sense one does obtain a weak solution. Concerning the first step of the above strategy we simply need to add to \cite{DR2} convective stability results for blocks with unstable endstates but no other kind of instability. 

With this in mind, we first observe that such a portion cannot possess a characteristic point because the closest characteristic point to the unstable infinity would then be an unstable characteristic point since the sign of $g'$ alternates along consecutive zeros of $g$. Therefore there are only two kinds of portions of wave profiles to consider: monostable non characteristic fronts and unstable constants on half-lines with characteristics leaving the half-line domain. Note that in contrast with \cite{DR1,DR2} we do need to consider some half-line problems, precisely because unstable constants are not convectively stable as whole line solutions. We stress that the condition that characteristics leave the half-line domain is consistent with entropy conditions (previously used to check that characteristics are outgoing from the unstable infinity), and that it also ensures that the corresponding initial boundary value problem may be solved without boundary conditions. 

Though all the kinds of refinements and extensions of \cite{DR1,DR2} --- extra discontinuities in perturbations, anisotropic smallness on gradients, higher regularity control, \emph{etc.} --- could also be obtained here, we restrain from providing such results since this would come with little new insights. The only exception we made is for the extension to multidimensional scalar equations. Our main motivation to deal explicitly with the multidimensional case, in Section~\ref{s:multiD}, is that the multidimensional analysis of \cite{DR1,DR2} is restricted to perturbations supported away from characteristic points and we take the present opportunity to extend it to general perturbations. The associated dynamics turns out to be significantly richer and, though it is not related to convective stability, we regard the identification of the large-time impacts arising from initial characteristic perturbations as a very important part of the present paper. We point out that the resolution of the latter does involve spatial weights in a way that may be thought as dual\footnote{In the sense that we estimate there (parts of) the solution with growing weights but measure initial data in weightless norms instead of estimating as here the solution with weightless norms while measuring initial data with decaying weights.} to the convective way.

As for functional space choices, as in \cite{DR1,DR2}, we work with the minimal amount of regularity enabling the propagation of further regularity, the --- here spatially weighted --- piecewise $W^{1,\infty}$ topology. Yet we restrict to data in the closure of smooth functions for this topology, the space, denoted as\footnote{For Bounded and Uniformly Continuous.} $BUC^1$, of functions that are piecewise $\cC^1$ with derivatives uniformly continuous and bounded (with weights). Note that by using $L^\infty$-based spaces --- as opposed to $L^p$ based spaces, $1\leq p<\infty$ --- one imposes no extra localization besides the one explicitly encoded by weights.

As for spatial weights, we mostly use exponential weights requiring solutions to decay as $\exp(-\kappa\,|x|)$ for some $\kappa>0$. The precise count is as follows. On the associated $L^\infty$-based spaces, the essential spectrum contains
\begin{align*}
\left\{\ -(f'(u_\infty)-\sigma)\iD\xi+g'(u_{\infty})-\kappa\,|f'(u_\infty)-\sigma|\ \middle|\ \xi\in\R\ \right\}\\
\vspace{5em}\,=\,\left\{\ \lambda\in\C\ \middle|\ \Re(\lambda)=g'(u_{\infty})-\kappa\,|f'(u_\infty)-\sigma|\ \right\}
\end{align*}
and this vertical line serves as one of the boundaries of the essential spectrum. Indeed, as we prove, there is a critical spatial decay growth 
\[
\kappa_{u_\infty}:=\frac{g'(u_{\infty})}{|f'(u_\infty)-\sigma|}
\]
such that when $\kappa>\kappa_{u_\infty}$ one obtains a spectral gap\footnote{Provided there is no other obstruction from another endstate/characteristic point/discontinuity point.} and exponential decay in time at the linear and nonlinear levels.

We also consider the critical case $\kappa=\kappa_{u_\infty}$ where the imaginary axis is included in the essential spectrum. For smooth waves, using this critical exponential weight alone, we prove bounded stability, with no decay rates. Note that the absence of decay rates is consistent with the fact that at the linearized level there is a fundamental obstruction to time decay without loss in the absence of a spectral gap; see the detailed discussion of the Datko-Pazy theorem in \cite{R_Roscoff}. Moreover $\exp(-\kappa_{u_\infty}|x|)$ is precisely the decay rate of the derivative of the profile $\uu$ (if not constant) thus the corresponding functional space allows for any bounded phase shift, including oscillator ones, leaving little hope for any decay to a uniform translate of the background wave; see the related analysis in \cite{Collet-Eckmann}. We prove that these two obstacles are the only ones. Explicitly, if $\rho$ is a positive nonincreasing function converging to zero at $+\infty$ not faster than exponentially\footnote{In a sense made precise in Definition~\ref{def:exp}.} and the perturbation decays as $\rho(|x|)\,\exp(-\kappa_{u_\infty}|x|)$, one may conclude time decay at rate $\rho(\theta\,t)$ for some $\theta>0$ with spatial decay $\exp(-\kappa_{u_\infty}|x|)$.

For discontinuous waves, we obtain almost identical results --- including bounded stability results with critical spatial decay $\exp(-\kappa_{u_\infty}|x|)$ --- provided that $[g(\uu)]_d/[\uu]_d<0$ at any discontinuity point $d$ of the background profile $\uu$ terminating an infinite portion with an unstable endstate, or equivalently in the generic case of Section~\ref{s:general} that no discontinuity point separates two constant parts of $\uu$ with one of those corresponding to an unstable endstate. Roughly speaking, the foregoing sign condition on the jump ratio ensures that the evolution of discontinuity locations obeys equations similar to affine equations with exponentially damping linear part and constant forcing by shape perturbations. This is sufficient to guarantee a uniform bounded control on discontinuity locations. In contrast, for phases encoding locations of discontinuities separating two constant parts of $\uu$, the linear damping is missing and this control is obtained by integrating in time bounds for their time derivatives that are exactly the same as the time-decay bounds for the shape perturbation. As a consequence, for waves exhibiting discontinuities separating two constant parts with at least one unstable endstate $u_\infty$, we prove stability results from data localized as $\rho(|x|)\,\exp(-\kappa_{u_\infty}|x|)$ under the extra assumption that $\rho$ is integrable since this spatial integrability does yield the required time integrability.

Actually another kind of departure from the smooth case occurs for discontinuous waves with at least one discontinuity point separating two non-constant portions. Such waves are not isolated as traveling waves (even if spatial translates are identified) and the asymptotic dynamics may lead in large time to \emph{another} nearby wave. 

Though we shall not elaborate corresponding details, we mention that an immediate application of our strategy also covers some degenerate cases when $g'(u_\infty)=0$. In this case the good direction of propagation of characteristics is not guaranteed and one must assume as a further assumption a good sign on $f'(u_\infty)-\sigma$. With this in hands, one may extend to this case all asymptotic stability results with integrable time decay rates, that is, either exponential in time with exponential spatial weights, or with rate $\rho(\theta\,t)$ for some $\theta>0$ starting from data localized as $\rho(|x|)$ where $\rho$ is sub-exponential and integrable. However our techniques cannot obtain for this case a bounded stability result from data with no localization.

Coming back to non degenerate waves, we point out that another interesting feature is that by increasing localization requirement on perturbations one may push arbitrary far from the imaginary axis the part of the essential spectrum directly due to an unstable endstate, so that in some sense no time-decay restriction arises from such endstates. In contrast, near a stable endstate the characteristics are incoming to the related infinity so that introducing spatial weights encoding extra localization would \emph{reduce} time-decay rates. As a consequence such stable endstates do come with maximal decay rates. Likewise, because of the sign of nearby characteristics, stable discontinuities bring obstruction to the increase of decay rates whereas stable critical point do not. Indeed for the latter improving the time decay is achieved by increasing the order of the initial vanishing near critical points.

\subsection{Perspectives  and outline}

After \cite{DR1,DR2} and the present contribution, we regard the study of \emph{degenerate} waves and the stability analysis of \emph{compound} patterns such as multi-traveling waves, in particular multi-fronts sometimes called terraces, as the main open questions concerning the stability of coherent structures of scalar balance laws. 

The development of a general stability theory for waves of \emph{general systems} of hyperbolic equations would obviously be very welcome. Yet we would like to point out that one cannot expect a simple spectral stability classification for general systems since already for constant states such a simple classification is impossible. However we do expect that a general Lyapunov theory, converting spectral stability into nonlinear stability, could be obtained for large classes of traveling waves. We refer the reader to \cite{FR1,FR2} for some first steps in this direction.

To conclude we would like to point out as another natural class of open questions the study of coherent structures of parabolic regularizations of balance laws approximating hyperbolic traveling waves in the \emph{vanishing viscosity limit}. The notion of solution, entropy solutions, that we use here coincides with the one of solutions that could be obtained as local-in-time limits of solutions to such regularized equations. It is therefore natural to expect that in a uniform neighborhood of stable traveling waves this vanishing viscosity approximation is indeed global in time. The apparition of boundary layers at discontinuity points makes it a dreadful multi-scale problem. The case of non-characteristic waves of scalar balance laws have been nevertheless analyzed in \cite{Blochas-Rodrigues,Blochas}.

\bigskip

In the next section we first consider separately non-characteristic monostable smooth fronts. In Section~\ref{s:general} we prove convective stability of suitable constants on half-lines and complete the stability classification. At last, in Secton~\ref{s:multiD} we show how to extend the results to scalar balance laws posed in higher dimensions and, when doing so, we complete the weightless analysis of \cite{DR1,DR2} so as to allow general perturbations of characteristic points.

\medskip

\noindent {\bf Acknowledgment:} L.M.R. would like to warmly thank Paul Blochas for enlightening discussions about critical decay cases during his preparation of \cite{Blochas} and the finalization of the present piece of work. L.G. expresses his gratitude to Université de Rennes for its hospitality during part of the preparation of the present contribution. 

\section{Smooth monostable fronts}
\label{s:fronts}

We begin our analysis with the consideration of non-characteristic monostable smooth fronts.

\subsection{Notation and results}

To reduce the number of cases we observe that for any choice of signs $\epsilon$ and $\epsilon'$ in $\{-1,1\}$, any smooth solution $u$ to \eqref{e:original-eq} generates a solution $(t,x)\mapsto \epsilon\,u(t,\epsilon'\,x)$ to an equation of the same form with flux $\epsilon\,\epsilon'\,f(\epsilon\,\cdot)$ and source $\epsilon\,g(\epsilon\,\cdot)$. Such a transformation sends a monostable smooth front of the original equation to one of the transformed equation. 

Thanks to this, without loss of generality we may focus on fronts with endstates\footnote{Note that we have dropped the $\infty$ sign for the sake of notational brevity.} $u_-$ and $u_+$ respectively at $-\infty$ and $+\infty$, such that $g'(u_-)<0<g'(u_+)$ and $u_-<u_+$.

As a preliminary we discuss the existence of the relevant fronts. Fix $(u_-,u_+)$ such that $g'(u_-)<0<g'(u_+)$, $u_-<u_+$ and there is no zero of $g$ in $(u_-,u_+)$. Then for any $\sigma>\max_{[u_-,u_+]}(f')$, one obtains a smooth traveling wave of speed $\sigma$ and profile $\uu$ by solving
\begin{align*}
\uu(0)&\,=\,\frac{u_-+u_+}{2}\,,&
\uu'&\,=\,\frac{g(\uu)}{f'(\uu)-\sigma}\,.
\end{align*}
All the other fronts with the same endstates (in the same order) and speed are associated with profiles $\uu(\cdot+x_0)$, $x_0\in\R$.

We now introduce the relevant functional spaces and norms. We recall that $BUC^1(\R)$ denotes the space of functions that are $\cC^1$, bounded and uniformly continuous, with first derivative bounded and uniformly continuous. For any $\kappa>0$, we introduce exponential weights $\Omega_\kappa$, $\Omega_\kappa^+$ defined by
\begin{align*}
\Omega_{\kappa}(x)&\,:=\,\exp(-\kappa\,x)\,,&
\Omega_{\kappa}^+(x)&\,:=\,\Omega_{\kappa}(x_+)\,,&
\end{align*}
where $(\cdot)_+$ denotes the positive part, that is, $a_+=\max(\{a,0\})$. We then set 
\[
\|v\|_{W^{1,\infty}_{\kappa,+}}
:=\|(\Omega_\kappa^+)^{-1}\,(v,\,\d_x v)\|_{L^\infty(\R)}\,.
\]
We recall that the critical spatial growth is 
\[
\kappa^+:=\frac{g'(u_+)}{|f'(u_+)-\sigma|}\,.
\]
We also introduce the corresponding time-decay exponential rates
\[
\omega_\kappa:=\min(\{|g'(u_-)|,(\kappa-\kappa^+)\,|f'(u_+)-\sigma|\})\,.
\]
Note that $\omega_\kappa>0$ when $\kappa>\kappa^+$.

Moreover we give a precise definition of what we call sub-exponential weight.

\bd\label{def:exp}
A sub-exponential weight is a positive function $\rho:\,\R_+\to\R_+^*$, non increasing, converging to $0$ at $+\infty$ and such that
\begin{enumerate}
\item there exist positive $\omega$ and $C$ such that for any $t\geq0$
\[
C\,\exp(-\omega\,t)\leq \rho(t)\,;
\]
\item there exist positive $\omega$ and $C$ such that for any $t\geq0$
\[
\int_0^t \exp(-\omega\,(t-s))\,\rho(s)\,\dD s\leq C\,\rho(t)\,.
\]
\end{enumerate} 
\ed

A simple way to satisfy both inequality constraints is to enforce that there exist positive $C'$ and $\omega'$ such that for any $0\leq s\leq t$
\[
\frac{\rho(s)}{\rho(t)}\leq C'\, \exp(\omega'(t-s))\,.
\]
In particular the class of sub-exponential weights includes decreasing exponentials $(\Omega_{\kappa})_{|\R_+}$, $\kappa>0$, but also algebraic weights such as
\[
\R_+\to\R_+^*\,,\quad x\mapsto \frac{1}{(1+x)^r}\,,\qquad\textrm{where}\quad r>0\,.
\]

With any such weight $\rho$ we associate $\rho^+:\R\to\R_+^*$ defined by $\rho^+(x)=\rho(x_+)$, and set 
\[
\|v\|_{W^{1,\infty}_{+,\rho}}
:=\|(\rho^+\,\Omega_{\kappa^+}^+)^{-1}\,(v,\,\d_x v)\|_{L^\infty(\R)}\,.
\]

\medskip

The goal of the present section is to prove the following statements.

\bt\label{th:fronts-exp}
Consider a smooth non-characteristic traveling wave of profile $\uu$ and speed $\sigma$ with endstates $u_-$ and $u_+$ respectively at $-\infty$ and $+\infty$ such that $g'(u_-)<0<g'(u_+)$.\\
Let $\kappa\geq\kappa^+$ where $\kappa^+:=g'(u_+)/|f'(u_+)-\sigma|$.\\
There exist positive $C$ and $\delta_0$ such that if $u_0\in BUC^1(\R)$ with $(u_0-\uu)/\Omega_\kappa^+\in BUC^0(\R)$,  $\d_x(u_0-\uu)/\Omega_\kappa^+\in BUC^0(\R)$ and
\[
E_{0,\kappa}:=\|u_0-\uu\|_{W^{1,\infty}_{\kappa,+}}
\leq \delta_0
\]
then the solution $u$ to \eqref{e:original-eq} starting from $u(0,\cdot)=u_0$ satisfies for any $t\geq0$ 
\begin{align*}
\|u(t,\cdot+\sigma\,t)-\uu\|_{W^{1,\infty}_{\kappa,+}}
\leq\,C\,E_{0,\kappa}\,
\exp(-t\,\omega_\kappa)
\end{align*}
where 
\[
\omega_\kappa:=\min(\{|g'(u_-)|,(\kappa-\kappa^+)\,|f'(u_+)-\sigma|\})\,.
\]
\et

Note that the above decay rates match exactly the spectral gaps of the linearized operator on the corresponding weighted space.

\bt\label{th:fronts}
Consider a smooth non-characteristic traveling wave of profile $\uu$ and speed $\sigma$ with endstates $u_-$ and $u_+$ respectively at $-\infty$ and $+\infty$ such that $g'(u_-)<0<g'(u_+)$. Set $\kappa^+:=g'(u_+)/|f'(u_+)-\sigma|$.\\
There exists a positive $\delta_0$ and for any sub-exponential weight $\rho$ there exist positive $C_\rho$ and $\theta_\rho$ such that if $u_0\in BUC^1(\R)$ with $(u_0-\uu)/\Omega_{\kappa^+}^+\in BUC^0(\R)$,  $\d_x(u_0-\uu)/\Omega_{\kappa^+}^+\in BUC^0(\R)$ and
\[
E_{0,\kappa^+}:=\|u_0-\uu\|_{W^{1,\infty}_{\kappa^+,+}}
\leq \delta_0
\]
then the solution $u$ to \eqref{e:original-eq} starting from $u(0,\cdot)=u_0$ satisfies 
\begin{enumerate}
\item if 
\be\label{as:extra}
\lim_{x\to\infty}\exp(\kappa^+x)\,(u_0-\uu,\d_x(u_0-\uu))(x)\,=\,(0,0)
\ee
then 
\[
\lim_{t\to\infty}\|u(t,\cdot+\sigma\,t)-\uu\|_{W^{1,\infty}_{\kappa^+,+}}\,=\,0\,;
\]
\item if $\rho$ is a sub-exponential weight such that 
\[
E_0^{[\rho]}:=\|u_0-\uu\|_{W^{1,\infty}_{+,\rho}}<+\infty
\]
then for any $t\geq0$
\begin{align*}
\|u(t,\cdot+\sigma\,t)-\uu\|_{W^{1,\infty}_{\kappa,+}}
\,\leq\,C_\rho\,E_0^{[\rho]}\,\rho(\,\theta_\rho\,t)\,.
\end{align*}
\end{enumerate}
\et

Note that when choosing as a sub-exponential weight an algebraic weight, the introduction of $\theta_\rho$ merely changes the constant in factor not the actual decay rate.

\br\textup{
An important feature of the decay result associated with sub-exponential weights is that the smallness assumption is required on $E_{0,\kappa^+}$ independently of $\rho$, but not on $E_0^{[\rho]}$ nor in a $\rho$-dependent way. This is the keystone of the proof of the rate-less decay result, which is obtained by designing a sub-exponential $\rho$ such that $E_0^{[\rho]}$ is finite but potentially very large and/or with little control on $\rho$. If we were assuming smallness on $E_0^{[\rho]}$ then the second inequality in the definition of sub-exponentiality could be relaxed to 
\[
\int_0^t \exp(-\omega\,(t-s))\,\rho(s)^2\,\dD s\leq C\,\rho(t)\,.
\]
Likewise in Theorem~\ref{th:fronts-exp} we may obtain the decay associated with $\kappa$ assuming the smallness $E_{0,\kappa^+}$ only (instead of $E_{0,\kappa}$) provided that $\kappa^+\leq \kappa\leq 2\kappa^+$.}
\er

\br\label{rk:xlimit}\textup{
If 
\[
\lim_{x\to\infty}\exp(\kappa^+x)\,(u_0-\uu,\d_x(u_0-\uu))(x)\,=\,(a_\infty,b_\infty)
\]
then necessarily $b_\infty=-\kappa^+\,a_\infty$. Moreover in this case there exists $\psi_\infty$ such that  
\[
\lim_{x\to\infty}\exp(\kappa^+x)\,(u_0-\uu(\cdot-\psi_\infty),\d_x(u_0-\uu(\cdot-\psi_\infty)))(x)\,=\,(0,0)
\]
so that, by using the invariance by spatial translation of \eqref{e:original-eq}, one may also deduce that
\[
\lim_{t\to\infty}\|u(t,\cdot+\sigma\,t)-\uu(\cdot-\psi_\infty)\|_{W^{1,\infty}_{\kappa^+,+}}\,=\,0\,.
\]
To be somewhat more explicit, setting $\alpha^+:=\lim_{x\to\infty} \exp(\kappa^+x)\,\uu'(x)$, one derives that
\[
\psi_\infty\,=\,\frac{1}{\kappa^+}\ln\left(1-\frac{\kappa^+}{\alpha^+}a_\infty\right)\,.
\]
}
\er

To motivate the linear analysis, we sketch the first steps of the proofs of Theorems~\ref{th:fronts-exp} and~\ref{th:fronts}. The goal is to bound $v(t,\cdot)=u(t,\cdot+\sigma\,t)-\uu$. This is done through the Duhamel-type formula
\begin{equation}\label{eq:Duhamel-fronts}
v(t,\cdot)\,=\,S_{f'(\uu+v)-\sigma}(0,t)(u_0-\uu)
+\int_0^tS_{f'(\uu+v)-\sigma}(s,t)\,N(v(s,\cdot))\,\dD s\,,
\end{equation}
where 
\begin{align*}
N(w)&:=g(\uu+w)-g(\uu)-g'(\uu)\,w
-\uu'\left(f'(\uu+w)-f'(\uu)-f''(\uu)\,w\right)\\
&\qquad-\frac{\uu''}{\uu'}(f'(\uu+w)-f'(\uu))\,w\,,
\end{align*}
and $S_a(t_1,t_2)$ denotes the evolution operator from time $t_1$ to time $t_2$ for the linear time-dependent equation
\[
\d_tw+a\,\d_xw\,=\,a\,\frac{\uu''}{\uu'}\,w\,.
\]
For the necessary background on evolution systems we refer to \cite[Chapter~5]{Pazy-83}. To gain in brevity, we skip below some of the immaterial details in the application of the corresponding theory and simply refer to \cite{DR1,DR2} for the omitted details.

\br\textup{
The motivation for the form of the zeroth order terms of the above linear equation is that the equation is then equivalently written as
\[
\d_t\left(\frac{w}{\uu'}\right)+a\,\d_x\left(\frac{w}{\uu'}\right)\,=\,0\,,
\] 
this form being more convenient when seeking for a modification of the first-order derivative of $w$ satisfying an equation in closed form. Other choices would typically induce a spurious loss of a $(1+t)$-factor (when the growth is saturated, that is, when $\omega_\kappa\geq |g'(u_-)|$).
}
\er

\subsection{Linear estimates}

At the linear level the goal is thus to bound $S_a(t_1,t_2)$ when $a$ is sufficiently close to $f'(\uu)-\sigma$. To prepare the proof of Theorem~\ref{th:fronts-exp} when $\kappa>\kappa^+$, the following proposition is sufficient.

\bpr\label{p:exp}
Under the assumptions of Theorem~\ref{th:fronts-exp}, for any $\kappa\geq\kappa^+$ there exist positive $\delta$ and $K$ such that if $T>0$ is such that
\[
\|a-(f'(\uu)-\sigma)\|_{L^\infty(0,T;W^{1,\infty}_{\kappa^+,+})}\leq \delta
\]
then for any $0\leq s\leq t\leq T$
\begin{align*}
\|S_a(s,t)w\|_{W^{1,\infty}_{\kappa,+}}
&\leq K\,\|w\|_{W^{1,\infty}_{\kappa,+}}
\,\eD^{-(t-s)\,\omega_\kappa}
\,\exp\left(K\,\int_s^t\|a(\tau,\cdot)-(f'(\uu)-\sigma)\|_{W^{1,\infty}_{\kappa^+,+}}\dD\tau\right)\,.
\end{align*}
\epr

\begin{proof}
We first prove a bound for $\|(\Omega_{\kappa}^+)^{-1}(S_a(s,t)w)\|_{L^\infty(\R)}$ in terms of $\|(\Omega_{\kappa}^+)^{-1}w\|_{L^\infty(\R)}$. As in \cite[Section~2.2]{DR1} and \cite[Section~4C]{DR2}, this follows from the abstract theory of evolution system provided we prove that there exist positive constants $K$ and $\delta$ and a weight $\alpha$ bounded from above and below by a positive multiple of $\Omega_\kappa^+$ such that for any $a_0$ (time-independent) and any $\lambda \in\C$ such that 
\begin{align*}
\|a_0-(f'(\uu)-\sigma)\|_{W^{1,\infty}_{\kappa^+,+}}&\leq \delta\,,&
\Re(\lambda)&> -\omega_\kappa+K\,\|a_0-(f'(\uu)-\sigma)\|_{W^{1,\infty}_{\kappa^+,+}}
\end{align*}
there holds that $\lambda$ does not belong to the spectrum of 
\[
L_{a_0}:=-a_0\,\d_x
+a_0\,\frac{\uu''}{\uu'}
\]
and
\[
\|\alpha^{-1}(\lambda\,\I-L_{a_0})^{-1}(w)\|_{L^\infty(\R)}\,
\leq\,
\frac{\|\alpha^{-1}\,w\|_{L^\infty(\R)}}{\Re(\lambda)+\omega_\kappa-K\,\|a_0-(f'(\uu)-\sigma)\|_{W^{1,\infty}_{\kappa^+,+}}}\,.
\]

To achieve such a goal, we seek $\alpha$ under the form $\alpha=\exp(\phi)$ for a Lipschitz $\phi$ so that
\[
\alpha^{-1}\,L_{a_0}\,w\,=\,
\left(-a_0\,\d_x-a_0\,\phi'
+a_0\,\frac{\uu''}{\uu'}
\right)(\alpha^{-1}w)\,.
\]
It is then convenient to define $\phi$ by 
\[
\phi'=-\frac{(\omega_\kappa+g'(\uu)-\uu'f''(\uu))_+}{|f'(\uu)-\sigma|}
\]
so that
\[
-a_0\,\phi'
+a_0\,\frac{\uu''}{\uu'}
\leq -\omega_\kappa
+|a_0-(f'(\uu)-\sigma)|\,
\left(\frac{(\omega_\kappa+g'(\uu)-\uu'f''(\uu))_+}{|f'(\uu)-\sigma|}
+\frac{|\uu''|}{|\uu'|}
\right)\,,
\]
which is sufficient to conclude by a direct integration of the resolvent equation.

We now move to bound $\|(\Omega_{\kappa}^+)^{-1}\d_x(S_a(s,t)w)\|_{L^\infty(\R)}$ in terms of $\|w\|_{W^{1,\infty}_{\kappa,+}}$. To do so we introduce $\tilde{v}$ defined by
\[
\tilde{v}(t,\cdot):=\d_x(S_a(s,t)w)-\frac{\uu''}{\uu'}S_a(s,t)w
\]
and observe that
\[
\d_t\tilde{v}+a\,\d_x\tilde{v}
\,=\,\left(a\frac{\uu''}{\uu'}-\d_xa\right)\tilde{v}\,.
\]
The arguments used above provide the same type of bound for the evolution operators of the latter equation as for $S_a$. This is sufficient to conclude the proof of the proposition.
%
\end{proof}

To complete the proof of Theorem~\ref{th:fronts-exp} when $\kappa=\kappa^+$, we need linear bounds that do not lose the integral factor depending on $a-(f'(\uu)-\sigma)$. Those are provided by the following proposition.

\bpr\label{p:critical}
Under the assumptions of Theorem~\ref{th:fronts-exp}, there exist positive $\delta$ and $K$ such that if $T>0$ is such that
\[
\|a-(f'(\uu)-\sigma)\|_{L^\infty(0,T;W^{1,\infty}_{\kappa^+,+})}\leq \delta
\]
then for any $0\leq s\leq t\leq T$
\begin{align*}
\|S_a(s,t)w\|_{W^{1,\infty}_{\kappa^+,+}}
&\leq K\,\|w\|_{W^{1,\infty}_{\kappa^+,+}}\,.
\end{align*}
\epr

\begin{proof}
The strategy of proof is similar to the one used for Proposition~\ref{p:exp}. The first departure is that when proving resolvent bounds to estimate $\|(\Omega_{\kappa}^+)^{-1}(S_a(s,t)w)\|_{L^\infty(\R)}$, we define $\phi$ through
\begin{align*}
\phi'&=-\left(-\frac{\uu''}{\uu'}\right)_+\,,&
\textrm{so that }&&\phi'-\frac{\uu''}{\uu'}&\leq0\,.
\end{align*}
 
A similar change must be made in the part of the proof bounding $\|(\Omega_{\kappa}^+)^{-1}\d_x(S_a(s,t)w)\|_{L^\infty(\R)}$. The linear equation to consider is still
\[
\d_t\tilde{v}+a\,\d_x\tilde{v}
\,=\,\left(a\frac{\uu''}{\uu'}-\d_xa\right)\tilde{v}\,.
\]
But its study requires an extra argument compared to the above. When deriving the resolvent bound with a frozen $a_0$, we first observe that there exists a positive $K'$ such that
\[
\left|\frac{\d_xa_0}{a_0}-\frac{f''(\uu)}{f'(\uu)-\sigma}\uu''\right|\,\leq\,K'\,\delta\,\Omega_{\kappa^+}
\]
uniform in $a_0$ when $a_0$ is $\delta$-close to $f'(\uu)-\sigma$. This allows to choose, as a convenient weight for the resolvent bound, $\alpha=\exp(\phi)$ with $\phi$ defined by
\begin{align*}
\phi'&=-\left(-\frac{\uu''}{\uu'}+\frac{f''(\uu)}{f'(\uu)-\sigma}\uu''
+K'\,\delta\,\Omega_{\kappa^+}\right)_+\,,&
\textrm{so that }&&
-a_0\,\left(\phi'-\frac{\uu''}{\uu'}\right)+\d_xa_0&\leq0\,,
\end{align*}
provided that $K'\delta<|g'(u_-)|$.
\end{proof}


We now turn to the other linear bounds used in the proof of Theorem~\ref{th:fronts}.

\bpr\label{p:rho}
Under the assumptions of Theorem~\ref{th:fronts}, there exists a positive $\delta$ and for any sub-exponential weight $\rho$, there exist positive $\eta_\rho$ and $K_\rho$ such that if $T>0$ is such that
\[
\|a-(f'(\uu)-\sigma)\|_{L^\infty(0,T;W^{1,\infty}_{\kappa^+,+})}\leq \delta
\]
then for any $0\leq s\leq t\leq T$ and any $0<\eta\leq \eta_\rho$
\begin{align*}
\|S_a(s,t)w\|_{W^{1,\infty}_{\kappa^+,+}}
&\leq K_\rho\,\rho(\eta\,(t-s))\,\|w\|_{W^{1,\infty}_{+,\rho}}
\,.
\end{align*}
\epr

\begin{proof}
We only prove a bound on $\|(\Omega_{\kappa}^+)^{-1}(S_a(s,t)w)\|_{L^\infty(\R)}$ in terms of $\|(\rho^+\,\Omega_{\kappa}^+)^{-1}w\|_{L^\infty(\R)}$. The necessary extension is then carried out exactly as in the proof of Proposition~\ref{p:critical}. 

As an intermediate step, we prove that if $a-(f'(\uu)-\sigma)$ is sufficiently small then for any sufficiently large $x_0$ when $\eta$ is sufficiently small and $K'$ is sufficiently large
\[
\|(\Omega_{\kappa}^+)^{-1}(S_a(s,t)w)\|_{L^\infty((x_0,\infty))}
\leq K'\,\rho(\eta\,(t-s))\,\|(\rho^+\,\Omega_{\kappa}^+)^{-1}w\|_{L^\infty((x_0,\infty))}\,.
\]
To achieve this intermediate goal, we use the Duhamel formula
\begin{align*}
&(\Omega_{\kappa^+}^{-1}S_a(s,t)w)_{|(x_0,+\infty)}\,=\,\cT_a(s,t)(\Omega_{\kappa^+}^{-1}w)_{|(x_0,+\infty)}\\
&+\int_s^t \cT_a(s,t)\,\left((g'(\uu)-g'(u_+)-\uu'f''(\uu)+\kappa^+(a(\tau,\cdot)-(f'(u_+)-\sigma)))
\Omega_{\kappa^+}^{-1}\,S_a(s,\tau)w\,\right)_{|(x_0,+\infty)}\hspace{-1em}\dD \tau
\end{align*}
where $\cT_a(t_1,t_2)$ denotes the solution operator from time $t_1$ to $t_2$ associated with
\[
\d_t w+a\,\d_x w\,=\,0
\]
for functions of the spatial variable $x\in (x_0,+\infty)$. Let $\Phi_a(t_1,t_2,x)$ denote the solution at time $t_2$ starting from $x$ at time $t_1$ to the differential equation $y'(t)=a(t,y(t))$. Then, since $a$ is negative, denoting $\theta_{\cT}$ a positive lower bound on $|a|$ uniform in $a$ when $a$ is sufficiently close to $f'(\uu)-\sigma$, one derives that when $t_1\geq t_2$
\[
\Phi_a(t_1,t_2,x)\geq x+\theta_{\cT}\,(t_1-t_2)
\]
thus for any $x\geq x_0>0$
\begin{align*}
|\cT_a(s,t)(w)(x)|=|w(\Phi_a(t,s,x))|&\leq \rho(x+\theta_{\cT}\,(t-s))\|\rho^{-1}w\|_{L^\infty((x_0,+\infty))}\\
&\leq \rho(\theta_{\cT}\,(t-s))\|\rho^{-1}w\|_{L^\infty((x_0,+\infty))}
\end{align*}
so that for any $0<\theta\leq \theta_{\cT}$
\[
\|\cT_a(s,t)(w)\|_{L^\infty((x_0,+\infty))}
\leq \rho(\theta\,(t-s))\ \|\rho^{-1}w\|_{L^\infty((x_0,+\infty))}\,.
\]
Applying this both to $\rho$ and an exponential weight, we deduce through the above Duhamel formula that for some positive $\kappa'$, $K'$ and $\omega'$ and for any $0<\theta\leq \theta_{\cT}$
\begin{align*}
\|\Omega_{\kappa^+}^{-1}S_a(s,t)w\|_{L^\infty((x_0,+\infty))}\,
&\leq\,\rho(\theta\,(t-s))\|(\rho^+\Omega_{\kappa^+})^{-1}w\|_{L^\infty((x_0,+\infty))}\\
&\qquad+K'\,\eD^{-\kappa'x_0}\,\int_s^t \eD^{-\omega'(t-\tau)}\,\|\Omega_{\kappa^+}^{-1}\,S_a(s,\tau)w\|_{L^\infty((x_0,+\infty))}\dD \tau\,,
\end{align*}
which leads through a Gr\"onwall lemma\footnote{That is, a comparison with the solution $z$ to 
\begin{align*}
z(t)\,
&=\,\rho(\theta\,(t-s))\|(\rho^+\Omega_{\kappa^+})^{-1}w\|_{L^\infty((x_0,+\infty))}
+K'\,\eD^{-\kappa'x_0}\,\int_s^t \eD^{-\omega'(t-\tau)}\,z(\tau)\dD \tau\,.
\end{align*}} to 
\begin{align*}
&\|\Omega_{\kappa^+}^{-1}S_a(s,t)w\|_{L^\infty((x_0,+\infty))}\,\\
&
\leq\,\left(\rho(\theta\,(t-s))
+K'\,\eD^{-\kappa'x_0}\int_s^t
\eD^{-(t-\tau)\left(\omega'-K'\,\eD^{-\kappa'x_0}\right)}
\rho(\theta\,(\tau-s))\dD \tau\right)
\|(\rho^+\Omega_{\kappa^+})^{-1}w\|_{L^\infty((x_0,+\infty))}\,.
\end{align*}
This is sufficient to conclude the proof of the intermediate step by taking $x_0$ sufficiently large and decreasing the upper bound on $\theta$ so as to use the sub-exponential character of $\rho$.

We now conclude the proof of the proposition by choosing a cut-off function $\chi$, $0\leq \chi\leq 1$, taking the value $1$ on $(-\infty,x_0)$ and vanishing on a neighborhood of $+\infty$ so that for some positive constants $K''$ and $\omega''$ when $\theta$ is positive and sufficiently small
\begin{align*}
\|&\Omega_{\kappa^+}^{-1}S_a(s,t)w\|_{L^\infty(\R)}\\
&=\Big\|\Omega_{\kappa^+}^{-1}S_a\Big(s+\frac{t-s}{2},t\Big)
\left(\chi S_a\Big(s,s+\frac{t-s}{2}\Big)w+(1-\chi)S_a\Big(s,s+\frac{t-s}{2}\Big)w\right)\Big\|_{L^\infty(\R)}\\
&\leq\ K''\,\eD^{-\omega''\frac{t-s}{2}}\Big\|\Omega_{\kappa^+}^{-1}S_a\Big(s,s+\frac{t-s}{2}\Big)w\Big\|_{L^\infty(\R)}
+K''\,\Big\|\Omega_{\kappa^+}^{-1}S_a\Big(s,s+\frac{t-s}{2}\Big)w\Big\|_{L^\infty((x_0,\infty))}\\
&\leq\ (K'')^2\,\eD^{-\omega''\frac{t-s}{2}}\|\Omega_{\kappa^+}^{-1}w\|_{L^\infty(\R)}
+(K'')^2\rho\Big(\theta\,\frac{t-s}{2}\Big)\,\|(\rho^+\Omega_{\kappa^+})^{-1}w\|_{L^\infty((x_0,+\infty))}\,.
\end{align*}
This is sufficient to derive the claimed bound.
\end{proof}

\subsection{Nonlinear stability - Proofs of Theorems~\ref{th:fronts-exp} and~\ref{th:fronts}}

\subsubsection*{Proof of Theorem~\ref{th:fronts-exp}}

The proof of Theorem~\ref{th:fronts-exp} follows from Propositions~\ref{p:exp} and~\ref{p:critical} and Duhamel formula \eqref{eq:Duhamel-fronts} through a continuity argument. The detailed proof cannot be identical for cases $\omega_\kappa\geq |g'(u_-)|$ and $\kappa=\kappa^+$, the remaining cases $0<\omega_\kappa<|g'(u_-)|$ being however treatable by one variant or the other. 

We first provide details on the case $\kappa>\kappa^+$. In this case there exist positive constants $K$ and $\delta$ such that if $2K\,E_{0,\kappa}\leq \delta$ then, when for some $T>0$, for any $0\leq t\leq T$,
\begin{align}\label{continuity}
\|v(t,\cdot)\|_{W^{1,\infty}_{\kappa,+}}
\leq\,2\,K\,E_{0,\kappa}\,\exp(-t\,\omega_\kappa)
\end{align}
there holds for any $0\leq t\leq T$,
\begin{align*}
\|v(t,\cdot)\|_{W^{1,\infty}_{\kappa,+}}
&\leq\,K\,E_{0,\kappa}\,\eD^{-t\,\omega_\kappa}\,\exp\left(\frac{2K^2}{\omega_\kappa}E_{0,\kappa}\right)\\
&\qquad+4\,K^3\,E_{0,\kappa}^2\,\exp\left(\frac{2K^2}{\omega_\kappa}E_{0,\kappa}\right)
\int_0^t\eD^{-(t-\tau)\,\omega_\kappa}
\eD^{-2\,\tau\,\omega_\kappa}\,\dD\tau\\
&\leq K\,E_{0,\kappa}\,\eD^{-t\,\omega_\kappa}\,
\exp\left(\frac{2K^2}{\omega_\kappa}E_{0,\kappa}\right)
\,\left(1+\frac{4\,K^2\,E_{0,\kappa}}{\omega_\kappa}\right)\,.
\end{align*}
Restricting further the size of $E_{0,\kappa}$ to enforce 
\[
\exp\left(\frac{2K^2}{\omega_\kappa}E_{0,\kappa}\right)
\,\left(1+\frac{4\,K^2\,E_{0,\kappa}}{\omega_\kappa}\right)<2
\]
we deduce that the set of times $T$ such that \eqref{continuity} holds for $0\leq t\leq T$ is both closed (by continuity alone) and open (by continuity and the above). Since this set contains $0$ it is equal to $\R_+$. Hence the result.

We now give details on the case $\kappa=\kappa^+$. In this case there exist positive constants $K$ and $\delta$ such that if $2K\,E_{0,\kappa}\leq \delta$ then, when for some $T>0$, for any $0\leq t\leq T$,
\begin{align}\label{cont-critical}
\|v(t,\cdot)\|_{W^{1,\infty}_{\kappa^+,+}}
\leq\,2\,K\,E_{0,\kappa^+}
\end{align}
there holds for any $0\leq t\leq T$,
\begin{align*}
\|v(t,\cdot)\|_{W^{1,\infty}_{\kappa^+,+}}
&\leq\,K\,E_{0,\kappa^+}\,
+4\,K^3\,E_{0,\kappa^+}^2\,
\int_0^t\eD^{-(t-\tau)\,(\omega_{2\kappa^+}-2K^2\,E_{0,\kappa^+})}\,\dD\tau\\
&\leq K\,E_{0,\kappa^+}
\,\left(1+\frac{4\,K^2\,E_{0,\kappa^+}}{\omega_{2\kappa^+}-2K^2\,E_{0,\kappa^+}}\right)\,,
\end{align*}
provided that furthermore $2K^2\,E_{0,\kappa^+}<\omega_{2\kappa^+}$. The proof is then achieved by continuity under the extra condition 
\[
\frac{4\,K^2\,E_{0,\kappa^+}}{\omega_{2\kappa^+}-2K^2\,E_{0,\kappa}}<1\,.
\]
This achieves the proof of Theorem~\ref{th:fronts-exp}.\hfill\qed

\subsubsection*{Proof of Theorem~\ref{th:fronts}}

We begin the proof of Theorem~\ref{th:fronts} with the $\rho$-bound. By using the $\kappa=\kappa^+$ bounds from Theorem~\ref{th:fronts-exp} and Propositions~\ref{p:exp} and~\ref{p:rho} and Duhamel formula \eqref{eq:Duhamel-fronts}, one derives that there exist positive $\delta$, $K$ and $\omega'$ and for any sub-exponential $\rho$ there exist positive $\eta_\rho$ and $K_\rho$ such that if $E_{0,\kappa^+}\leq \delta$ then for any $t\geq 0$ and any $0<\rho\leq \eta_\rho$
\begin{align*}
\|v(t,\cdot)\|_{W^{1,\infty}_{\kappa^+,+}}
&\leq\,K_\rho\,\rho(\theta\,t)\,E_0^{[\rho]}\,
+K\,E_{0,\kappa^+}\,
\int_0^t\eD^{-(t-\tau)\,\omega'}
\|v(\tau,\cdot)\|_{W^{1,\infty}_{\kappa^+,+}}\,\dD\tau\,,
\end{align*}
which leads through a Gr\"onwall lemma to for any $t\geq 0$ and any $0<\rho\leq \eta_\rho$
\begin{align*}
\|v(t,\cdot)\|_{W^{1,\infty}_{\kappa^+,+}}
&\leq\,K_\rho\,E_0^{[\rho]}\,
\left(\,\rho(\theta\,t)
+K\,E_{0,\kappa^+}\,
\int_0^t\eD^{-(t-\tau)\,(\omega'-K\,E_{0,\kappa^+})}
\rho(\theta\,\tau)\,\dD\tau\right)\,.
\end{align*}
By imposing the further restriction $K\,E_{0,\kappa^+}<\omega'$ (independently of $\rho)$, one concludes the proof of the $\rho$-bound thanks to the sub-exponential character of $\rho$.

We now turn to the rate-less part of Theorem~\ref{th:fronts}. Assuming 
\[
\lim_{x\to\infty}\exp(\kappa^+x)\,(u_0-\uu,\d_x(u_0-\uu))(x)\,=\,0
\]
we shall build a sub-exponential weight $\rho$ such that $E_0^{[\rho]}<+\infty$. Let us first build recursively $(x_j)_{j\in\N}$ such that $x_0=0$ and for any $j$, $x_{j+1}\geq x_j+1$ and for any $x\geq x_{j+1}$, 
\[
|\exp(\kappa^+x)\,(u_0-\uu,\d_x(u_0-\uu))(x)|\leq \frac{1}{2^{j+1}}\,E_{0,\kappa^+}\,.
\] 
Then we define 
\[
\rho\,:=\,\sum_{j\in \N} \frac{1}{2^j}\,\One_{[t_j,t_{j+1})}
\]
where $\One_A$ denotes the characteristic function of a set $A$. By design $E_0^{[\rho]}\leq E_{0,\kappa^+}<+\infty$. We now check that $\rho$ is indeed sub-exponential. To do so we introduce the map $J:\R_+\to \N$ that associates with $x$ the single $j$ such that $x\in [x_j,x_{j+1})$. We observe that for any $x\leq y$, $J(y)-J(x)\leq y-x+2$. This implies that for any $x\leq y$
\[
\frac{\rho(x)}{\rho(y)}\,=\,2^{J(y)-J(x)}\,\leq\,4\,\exp((y-x)\ln(2))\,.
\]
As already pointed out this is sufficient to check that $\rho$ is sub-exponential. Applying the $\rho$-bound then concludes the proof.\hfill\qed

\section{General classification}
\label{s:general}

We now turn to general waves of \eqref{e:original-eq} that are non degenerate in the sense of Definition~\ref{def:non-degenerate}.

\subsection{Generic equations}

To provide a classification of stable waves as compact as possible, as in \cite{DR2} we introduce the following assumption.

\begin{assumption}\label{as:generic}
For any $(u,v)\in\R^2$, if
\begin{align*}
g(u)&=0\,,&g(v)&=0\,,&
g'(u)&\geq 0\,,&g'(v)&\geq 0\,,&
\textrm{and }&&f'(u)&=f'(v)\,,
\end{align*}
then $u=v$.
\end{assumption}

We refer to \cite{DR2} for a discussion of the precise sense in which Assumption~\ref{as:generic} (denoted Assumption~3.1 there) is generic.

Some comments on the notions of stability are necessary. As in the previous section we use a $BUC^1$ space with corresponding $W^{1,\infty}$-norm. However functions are considered on $\R\setminus D$ where $D$ is the set of discontinuity points of the background profile $\uu$ and uniform continuity is required separately on each connected component of $\R\setminus D$. We shall say that \emph{convective} stability holds if it stability holds when measured with a norm of the form 
\[
\|v\|_{W^{1,\infty}_\Omega(\R\setminus D)}\,:=\,\|\Omega^{-1}\,(v,\d_xv)\|_{L^\infty}
\]
where $\Omega$ is a positive weight, belonging to $BUC^0(\R\setminus D)$.

Moreover we measure proximity in a sense adapted to general traveling waves, coined as \emph{space-modulated} stability in \cite{JNRZ-periodic}. The reader is referred to \cite{R,R_Roscoff} for a detailed informal motivation of the notion, including some discussion about why it is necessary and sharp. It consists here in measuring the shape deformation from the solution at time $t$, $u(t,\cdot+\sigma\,t)$ to the background profile $\uu$  with $\delta_{sm}(u(t,\cdot+\sigma\,t),\uu)$ where
\[
\delta_{sm}(v_1,v_2):=\inf_{\Psi\ \cC^1\ \textrm{diffeomorphsim}}
(\|v_1\circ \Psi-v_2\|_{W^{1,\infty}_\Omega(\R\setminus D)}
+\|\d_x(\Psi-\Id_\R)\|_{W^{1,\infty}_\Omega(\R)})\,.
\]
When $\uu$ is constant, it is sufficient to measure proximity with $\|u(t,\cdot+\sigma\,t)-\uu\|_{W^{1,\infty}_\Omega(\R)}$, which corresponds to the classical notion of stability. When $D$ is reduced to a single point, it is sufficient to use the well-known notion of \emph{orbital} stability, which in the present case controls $\delta_{orb}(u(t,\cdot+\sigma\,t),\uu)$ where
\[
\delta_{orb}(v_1,v_2):=\inf_{\varphi\in\R}
\|v_1(\cdot+\varphi)-v_2\|_{W^{1,\infty}_\Omega(\R\setminus D)}\,.
\]
In the foregoing discussion and elsewhere our convention is that the norm of a function not belonging to the associated functional space is $+\infty$.

\bt\label{th:classification}
Under Assumption~\ref{as:generic}, for any non-degenerate wave of \eqref{e:original-eq} in the sense of Definition~\ref{def:non-degenerate} convective spectral stability and convective nonlinear stability coincide and they are equivalent to the wave being defined by some $(\uu,\sigma)$ of one of the following forms
\begin{enumerate}
\item $\uu$ is constant of value $u_\infty$ such that $g'(u_\infty)<0$.
\item $\uu$ is a smooth front with no characteristic point. In this case, $\uu$ is a monostable front.
\item $\uu$ is a smooth front with exactly one characteristic point and its characteristic value $u_\star$ satisfies $g'(u_\star)>0$. In this case, $\uu$ is a bistable front.
\item $\uu$ possesses exactly one or two discontinuity points, at most one characteristic point and $[g(\uu)]_{d}/[\uu]_{d}\leq 0$ at any of the discontinuity points $d$. 
\end{enumerate}
Moreover for convectively stable waves with no discontinuity between two non-constant portions convective nonlinear asymptotic stability holds whereas for other convectively stable waves any solution arising from a small perturbation converges to a possibly different nearby wave.
\et

We stress that the last case gathers many different shape profiles. Indeed for a single discontinuity point the smooth parts could be constant, non constant non characteristic or non constant with a single characteristic point. For two discontinuity points, the middle part is necessarily non constant with a single characteristic point whereas the two infinite parts are non characteristic but could be either constant or non constant.

Compared to the weightless statement \cite[Theorem~3.2]{DR2}, what we have gained by going to \emph{convective} stability, besides the relaxation of some sign constraints, is the possibility to include waves with a single discontinuity but no constant part and waves with two discontinuities that do not have two constant infinite parts. These are precisely the kinds of waves that are not isolated.

We shall not give details about the instability part of Theorem~\ref{th:classification} since it follows immediately from the proof of the corresponding one of \cite[Theorem~3.2]{DR2} and the fact that we have enforced in our definition of convective stability that involved weights $w$ should be bounded. The rest of the section is devoted to its stability part that is not already covered by \cite[Theorem~3.2]{DR2}.

\subsection{Convectively stable constant states on half lines}

We recall that our strategy of proof for stability consists in gluing together smooth parts and that the analysis of these smooth parts require only the consideration of the convective stability of two new kinds of relative equilibria: the non-characteristic monostable fronts analyzed in the previous section and constant states on half lines.

Suitable versions of Theorems~\ref{th:fronts-exp} and~\ref{th:fronts} also hold for the latter case, with essentially the same proof (in a slightly simpler version). Thus we merely explain what should be changed in their statements.

By using symmetries of the equation it is sufficient to consider the convective stability of $\uu\equiv u_\infty$ with $g(u_\infty)=0$ and $g'(u_\infty)>0$, as a solution to 
\be\label{eq:moving}
\partial_t w + \partial_x (f(w)-\sigma\,w) = g(w)
\ee
an equation for functions $w$, $w(t,x)\in\R$ depending on time variable $t$ and spatial variable $x\in (0,\infty)$. For this equation to be uniquely solvable near $u_\infty$ without boundary condition (and consistently with our later use), we assume $\sigma>f'(u_\infty)$.

At the basis of the well-posedness of the corresponding dynamics lies the fact that when $a_0$ is negatively valued and $\Re(\lambda)$ is sufficiently large, the resolvent problem
\[
\lambda\,v+a_0\,\d_xv\,=\,b_0\,v+F
\]
for functions over $(0,\infty)$ is uniquely solvable with solution given by
\[
v(x)\,=\,-\int_x^{+\infty}\eD^{\int_x^y\frac{\lambda-b_0}{a_0}} \frac{F(y)}{a_0(y)}\dD y\,.
\]

The statements of counterparts to Theorems~\ref{th:fronts-exp} and~\ref{th:fronts} are obtained by replacing function spaces over $\R$ with function spaces over $(0,\infty)$, $u_+$ with $u_\infty$, $u(t,\cdot+\sigma\,t)-\uu$ with $w(t,\cdot)-u_\infty$ and dropping anything related to $u_-$ so that for instance
\begin{align*}
\kappa^+&:=\frac{g'(u_\infty)}{|f'(u_\infty)-\sigma|}\,,&
\omega_\kappa&:=(\kappa-\kappa^+)\,|f'(u_\infty)-\sigma|\,.
\end{align*}

\br
In the introduction, we have motivated the need of the extra decay assumption \eqref{as:extra} in the critical case $\kappa=\kappa^+$ in order to deduce large-time decay by the fact that non-characteristic monostable fronts are not isolated relative equilibria in the topology of $\|\cdot\|_{W^{1,\infty}_{\kappa^+,+}}$. Indeed the family of its spatial translates contains elements that are arbitrarily close to the original wave and modulating those may create an initial perturbation that could be used to create wild dynamics. A similar phenomenon/obstruction occurs for the constant solution $\uu\equiv u_\infty$ on $(0,\infty)$. Yet, this time, the family of nearby equilibria is obtained by solving 
\begin{align*}
\uu_{u^{(0)}}(0)&\,=u^{(0)}\,,&
\uu_{u^{(0)}}'&\,=\,\frac{g(\uu_{u^{(0)}})}{f'(\uu_{u^{(0)}})-\sigma}
\end{align*}
with $u^{(0)}$ sufficiently close to $u_\infty$. Incidentally we point out that under some extra global assumptions\footnote{About the ordering of the vanishing of $g$ and $f'-\sigma$ near $u_\infty$.} $\uu_{u^{(0)}}$ could be obtained as tails of translates of a monostable front with $u_\infty$ as an unstable endstate at $+\infty$.
\er

\subsection{Reconstructing phases}

For the sake of brevity and readability, we prove only some representative parts of the stability content of Theorem~\ref{th:classification}. For the same reason, we only provide critical decay statements, that is, counterparts to Theorem~\ref{th:fronts} but not to Theorem~\ref{th:fronts-exp}.

We give two detailed results: one where the Rankine-Hugoniot equation does come with a linear decay mechanism  for the discontinuity location and another one where it does not and shape deformation must decay sufficiently fast to be integrable in time. The latter happens when $[g(\uu)]_{d}=0$ at the discontinuity point $d$, thus in particular when $d$ separates two constant parts. We begin with such an example.

\bt\label{th:Riemann}
Consider a non-degenerate traveling wave of profile $\uu$ and speed $\sigma$ such that $\uu_{|(-\infty,0))}\equiv u_-$ and  $\uu_{|(0,+\infty))}\equiv u_+$ with constant values $u_-$ such that $g'(u_-)<0<g'(u_+)$.\\
Set $\kappa^+:=g'(u_+)/|f'(u_+)-\sigma|$ and let $\rho$ be an integrable sub-exponential weight. Define weights $\Omega:\R^*\to\R$ and $\Omega_\rho:\R^*\to\R$ by
\begin{align*}
\Omega(x)&:=\exp(-\kappa^+\,x_+)\,,&
\Omega_\rho(x)&:=\rho(x_+)\exp(-\kappa^+\,x_+)\,. 
\end{align*}

There exist positive $C$, $\theta$ and $\delta_0$ such that if $u_0\in BUC^1(\R^*)$ with $(u_0-\uu)/\Omega_\rho\in BUC^0(\R^*)$,  $\d_x(u_0-\uu)/\Omega_\rho\in BUC^0(\R^*)$ and
\[
E_0^{[\rho]}:=\|u_0-\uu\|_{W^{1,\infty}_{\Omega_\rho}(\R^*)}
\leq \delta_0
\]
then there exists $\psi\in BUC^2(\R_+)$ such that the solution $u$ to \eqref{e:original-eq} starting from $u(0,\cdot)=u_0$ satisfies for any $t\geq0$ 
\begin{align*}
\|u(t,\cdot+\sigma\,t+\psi(t))-\uu\|_{W^{1,\infty}_\Omega(\R^*)}+|(\psi',\psi'')(t)|
\leq\,C\,E_0^{[\rho]}\,\rho(\theta\,t)
\end{align*}
and there exists $\psi_\infty\in\R$ such that for any $t\geq0$
\begin{align*}
|\psi_\infty|&\leq C\,E_0^{[\rho]}\,,&
|\psi(t)-\psi_\infty|
\leq\,C\,E_0^{[\rho]}\,\int_{\theta\,t}^\infty\rho\,.
\end{align*}
\et

\begin{proof}
The first step is to extend smooth parts and apply stability results for constants states. On one hand, one may extend $(u_0)|_{|(-\infty,0)}$ into $u_{0,\ell}$ defined on $\R$ so that 
\[
\|u_{0,\ell}-u_-\|_{W^{1,\infty}(\R)}\leq 2\,E_0^{[\rho]} 
\]
and, provided that $E_0^{[\rho]}$ is sufficiently small, deduce from \cite[Section~2]{DR1} a control on the solution $u_\ell$ to \eqref{e:original-eq} arising from $u_{0,\ell}$,
\[
\|u_\ell(t,\cdot+\sigma\,t)-u_-\|_{W^{1,\infty}(\R)}\leq K_\ell\,E_0^{[\rho]}\,\exp(g'(u_-)\,t)\,, 
\]
for some $K_\ell>0$ independent of $u_0$. On the other hand, for the sub-exponential weight $\tilde{\rho}:=\rho((\cdot-1)_+)$, one may extend $(u_0((\cdot-1)_+))_{|(1,\infty)}$ into $w_{0,r}$ defined on $\R_+$ such that 
\[
\|w_{0,r}-u_+\|_{W^{1,\infty}_{\Omega_{\tilde{\rho}}}((0,\infty))}\leq K_0\,E_0^{[\rho]} 
\]
for some $K_0$ independent of $u_0$ and apply the stability results of the last subsection so as to derive a control on the solution $w_r$ to \eqref{eq:moving} arising from $w_{0,r}$,
\[
\|w_r(t,\cdot)-u_+\|_{W^{1,\infty}(\R)}\leq K'\,E_0^{[\rho]}\,\tilde{\rho}(\theta'\,t)\,, 
\]
for some $K'>0$ and $\theta'>0$ independent of $u_0$. Introducing $u_r$ defined by $u_r(t,x):=w_r(t,x-\sigma\,t+1)$ when $x\geq \sigma\,t-1$, one deduces from the later
\[
\|u_r(t,\cdot+\sigma\,t)-u_+\|_{W^{1,\infty}_{\Omega}((-1,+\infty))}\leq K_r\,E_0^{[\rho]}\,\rho(\theta_r\,t)\,, 
\]
for some $K_r>0$ and $\theta_r>0$ independent of $u_0$. Here, to go from $\tilde{\rho}$ to $\rho$, we have used that if $0<\theta<\theta'$, then when $t$ is sufficiently large
\[
(\theta'\,t-1)_+=\theta't-1\geq \theta\,t\,.
\] 

The second step is to define $\psi$ through $\psi(0)=0$ and the Rankine-Hugoniot condition
\be\label{eq:RH}
\psi'(t)\,=\,F(u_r(t,\sigma\,t+\psi(t)),u_\ell(t,\sigma\,t+\psi(t))\,)
-F(u_+,u_-\,)\,,
\ee
where
\be\label{eq:F}
F(v_1,v_2):=\int_0^1\,f'(\tau\,v_1+(1-\tau)\,v_2)\,\dD \tau\,,
\ee
and represent $u$ as 
\be\label{eq:glue}
u(t,x)\,=\,
\begin{cases}
u_\ell(t,x)&\quad\textrm{if }x< \sigma\,t+\psi(t)\,,\\
u_r(t,x)&\quad\textrm{if }x>\sigma\,t+\psi(t)\,.
\end{cases}
\ee
At this stage we need to rule out the possibility that \eqref{eq:RH} could lead to a finite-time blow-up by $\psi$ reaching the value $-1$. This is where we use the integrability of $\rho$ since \eqref{eq:RH} readily implies that $|\psi'(t)|$ is bounded by a multiple of $\rho(\theta\,t)\,E_0^{[\rho]}$ for some $\theta>0$. From this observation the rest of the bounds follows in a straightforward way, with
\[
\psi_\infty:=\int_0^{+\infty}\psi'\,.
\]
\end{proof}

We now give an example with a linear dissipation at the discontinuity. Moreover we choose an example with a discontinuity point separating two non-constant parts to illustrate convergence to another nearby wave. 

As a preparation, we give some details about how nearby waves arise, specializing to the kind of waves considered below. Let $\uu$ be a wave profile of speed $\sigma$ with exactly one discontinuity point at $0$, separating two non-constant smooth parts and assume that $[g(\uu)]_0/[\uu]_0<0$. First we note that we may solve 
\begin{align*}
\uu_\ell(0)&=\uu(0^-)\,,&
\uu_\ell'&=\frac{g(\uu_\ell)}{f'(\uu_\ell)-\sigma}
\end{align*}
on an open interval containing $(-\infty,0]$ and that $(\uu_\ell)_{|(-\infty,0)}=\uu_{|(-\infty,0)}$. Likewise we may obtain $\uu_r$ defined on an open interval containing $[0,+\infty)$ and extending $\uu_{|(0,+\infty)}$. Then by the implicit function theorem for any sufficiently small $\varphi_\ell$ there exists a $\psi_r$ smoothly depending on $\varphi_\ell$ such that 
\[
\R^*\to\R\,,\qquad x\mapsto\begin{cases}
\uu_\ell(x+\psi_r-\varphi_\ell)&\quad \textrm{if }x<0\,,\\
\uu_r(x+\psi_r)&\quad \textrm{if }x>0\,,
\end{cases}
\]
defines a wave profile of speed $\sigma$. To be more explicit, the equation for $\psi_r$ is
\[
(f(\uu_r)-\sigma\,\uu_r)(\psi_r)\,=\,(f(\uu_\ell)-\sigma\,\uu_\ell)(\psi_r-\varphi_\ell)
\]
and the condition ensuring local invertibility is $(f(\uu_r)-\sigma\,\uu_r)'(0)\neq(f(\uu_\ell)-\sigma\,\uu_\ell)'(0)$ (which is equivalent to $[g(\uu)]_0\neq0$).

\bt\label{th:mixed}
Consider a non-degenerate traveling wave of profile $\uu$ and speed $\sigma$ such that $\uu_{|(-\infty,0))}$ is smooth with a single characteristic point, denoted $x_\star$, and  $\uu_{|(0,+\infty))}$ is smooth non constant with no characteristic point and $\uu$ jumps at $0$ with $[g(\uu)]_0/[\uu]_0<0$. In particular\footnote{Recall that the sign of $f'(\uu)-\sigma$ is known near a discontinuity point and that signs of $f'(\uu)-\sigma$ and $g'(\uu)$ are related near infinities.}, denoting as $u_-$ and $u_+$ the corresponding endstates, $g'(u_-)<0<g'(u_+)$.\\
Set $\kappa^+:=g'(u_+)/|f'(u_+)-\sigma|$. Define a weight $\Omega:\R^*\to\R$ and, with any sub-exponential weight $\rho$, associate a weight $\Omega_\rho:\R^*\to\R$ by
\begin{align*}
\Omega(x)&:=\exp(-\kappa^+\,x_+)\,,&
\Omega_\rho(x)&:=\rho(x_+)\exp(-\kappa^+\,x_+)\,.
\end{align*}

There exists a positive $\delta_0$ and for any sub-exponential weight $\rho$ there exist positive $C_\rho$ and $\theta_\rho$ such that if $u_0\in BUC^1(\R^*)$ with $(u_0-\uu)/\Omega\in BUC^0(\R^*)$,  $\d_x(u_0-\uu)/\Omega\in BUC^0(\R^*)$ and
\[
E_0:=\|u_0-\uu\|_{W^{1,\infty}_\Omega(\R^*)}
\leq \delta_0
\]
then there exists $\psi\in BUC^2(\R_+)$ such that the solution $u$ to \eqref{e:original-eq} starting from $u(0,\cdot)=u_0$ satisfies for any $t\geq0$ 
\begin{align*}
\|u(t,\cdot+\sigma\,t+\psi(t))-\uu\|_{W^{1,\infty}_\Omega(\R^*)}+|(\psi,\psi',\psi'')(t)|
\leq\,C\,E_0
\end{align*}
and, denoting $(\cU^{[y]})_y$ the smooth family of nearby profiles with speed $\sigma$ normalized by the fact that the discontinuity point of $\cU^{[y]}$ is at $0$ and its characteristic point is at $y$,
\begin{enumerate}
\item if 
\[
\lim_{x\to\infty}\exp(\kappa^+x)\,(u_0-\uu,\d_x(u_0-\uu))(x)\,=\,(0,0)
\]
then there exists $\psi_\star$ such that $|\psi_\star|\leq C\,E_0$ and
\[
\lim_{t\to\infty}\|u(t,\cdot+\sigma\,t+\psi(t))-\cU^{[x_\star+\psi_\star]}\|_{W^{1,\infty}_\Omega(\R^*)}\,=\,0\,,
\]
and there exists $\psi_\infty$ such that
\[
\lim_{t\to\infty}(\psi,\psi',\psi'')(t)\,=\,(\psi_\infty,0,0)\,;
\]
\item if $\rho$ is a sub-exponential weight such that 
\[
E_0^{[\rho]}:=\|u_0-\uu\|_{W^{1,\infty}_{\Omega_\rho}(\R^*)}<+\infty
\]
then for $(\psi_\star,\psi_\infty)$ as above and any $t\geq0$
\begin{align*}
\|u(t,\cdot+\sigma\,t+\psi(t))-\cU^{[x_\star+\psi_\star]}\|_{W^{1,\infty}_\Omega(\R^*)}
+|(\psi-\psi_\infty,\psi',\psi'')(t)|
\,\leq\,C_\rho\,E_0^{[\rho]}\,\rho(\,\theta_\rho\,t)\,.
\end{align*}
\end{enumerate}
\et

The content of Remark~\ref{rk:xlimit} is also applicable here up to obvious changes, including that a different limit at $+\infty$ yields convergence not necessarily to a translate of $\uu$ but to some nearby wave.

\begin{proof}
The first step is again the extension step. Now we also need to extend the background profile. First, let $I_\ell$ be a neighborhood of the closure of $\uu((-\infty,0))$ on which $f'-\sigma$ and $g$ vanish exactly once, extend $f_{|I_\ell}$ and $g_{|I_\ell}$ into $f_\ell$ and $g_\ell$ so that $\uu_{|(-\infty,0)}$ is the restriction of a non degenerate smooth wave profile $\uu_\ell$ of speed $\sigma$ for
\be\label{eq:left}
\d_tu+\d_x(f_\ell(u))=g_\ell(u)
\ee
with a single characteristic point and an endstate $u_{+,\ell}$ at $+\infty$ such that $g'(u_{+,\ell})=g'(u_-)$, and then extend $(u_0)_{|(-\infty,0)}$ into $u_{0,\ell}$ defined on $\R$ so that 
\[
\|u_{0,\ell}-\uu_\ell\|_{W^{1,\infty}(\R)}\leq 2\,E_0\,. 
\]
Applying the results of \cite[Section~4]{DR2} one receives that there exists $\psi_{char}\in\R$ so that the solution $u_\ell$ to \eqref{eq:left} arising from $u_{0,\ell}$, satisfies for any $t\geq0$
\[
\|u_\ell(t,\cdot+\sigma\,t+\psi_{char})-\uu_\ell\|_{W^{1,\infty}(\R)}\leq K_\ell\,E_0\,\exp(-\min(\{|g'(u_-)|,g'(u_\star)\})\,t)\,, 
\]
with $|\psi_{char}|\leq K_\ell E_0$, for some $K_\ell>0$ independent of $u_0$, where $u_\star$ denotes the characteristic value taken by $\uu_{|(-\infty,0)}$. Incidentally we point out that $\psi_{char}$ is determined by $u_{0,\ell}(x_\star+\psi_{char})=u_\star$.  
Likewise we obtain $I_r$, $f_r$, $g_r$, $\uu_r$, $u_{0,r}$, and apply the results of the preceding section to derive 
\[
\|u_r(t,\cdot+\sigma\,t)-\uu_r\|_{W^{1,\infty}_{\Omega}(\R)}\leq K_r\,E_0\,,
\]
$\uu_r$ being a smooth non-characteristic monostable front of 
\[
\d_tu+\d_x(f_r(u))=g_r(u)\,.
\]

In the second step we define $\psi$ by the Rankine-Hugoniot equation
\begin{align*}
\psi'(t)\,&=F(\uu_r(\psi(t)),\uu_\ell(\psi(t))\,)
-F(\uu_r(0),\uu_\ell(0)\,)\\
&\quad+\,F(u_r(t,\sigma\,t+\psi(t)),u_\ell(t,\sigma\,t+\psi(t))\,)
-F(\uu_r(\psi(t)),\uu_\ell(\psi(t))\,)\,,
\end{align*}
where $F$ is as in \eqref{eq:F} and represent $u$ as in \eqref{eq:glue}. Here the possible blow-up would occur when $\sigma\,t+\psi(t)$ leaves $u_\ell(t,\cdot)^{-1}(I_\ell)\cap u_r(t,\cdot)^{-1}(I_r)$ or $\psi(t)$ leaves $\uu_\ell^{-1}(I_\ell)\cap \uu_r^{-1}(I_r)$. To enforce \eqref{eq:glue} we also need to ensure that for any $x\leq \sigma\,t+\psi(t)$, $u_\ell(t,x)\in I_\ell$ and for any $x\geq \sigma\,t+\psi(t)$, $u_r(t,x)\in I_r$. All those could be obtained by maintaining $\psi(t)$ sufficiently small (and taking $E_0$ sufficiently small).  Note that 
\[
(F(\uu_r(\cdot),\uu_\ell(\cdot)\,))'(0)
\,=\,\frac{[(f(\uu))']_0[\uu]_0-[f(\uu)]_0[\uu']_0}{([\uu]_0)^2}
\,=\,\frac{[(f(\uu)-\sigma\uu)']_0}{[\uu]_0}
\,=\,\frac{[g(\uu)]_0}{[\uu]_0}<0
\]
so that there exist positive $\delta_{RH}$, $\omega_{RH}$ and $K_{RH}'$, independent of $u_0$, such that if $\max_{[0,t]}|\psi|\leq \delta_{RH}$ then 
\[
|\psi(t)|\leq K_{RH}'E_0\,\int_0^t \eD^{-\omega\,(t-s)} \dD s\,.
\]
Thus if $E_0$ is sufficiently small we deduce by a continuity argument that the above constraints are globally satisfied and for any $t\geq 0$, $|\psi(t)|\leq K_{RH}\,E_0$ for some $K_{RH}$ independent of $u_0$. The rest of the decayless part is then readily obtained. 

We now move to the second part of the statement. The $\rho$-less statement is deduced from the $\rho$ statement as in the previous section. Concerning the latter, by the analysis of the previous section, we obtain that 
\[
\|u_r(t,\cdot+\sigma\,t)-\uu_r\|_{W^{1,\infty}_{\Omega}(\R)}\leq K_\rho\,E_0^{[\rho]}\,\rho(\theta_\rho'\,t)\,,
\]
for some positive $K_\rho$ and $\theta_\rho$ independent of $u_0$. Then we apply the construction of nearby profiles preceding the theorem with $\varphi_\ell:=\psi_{char}$ so as to receive a convenient $\psi_r$ and set $\psi_\star:=\psi_{char}-\psi_r$, $\psi_\infty:=\psi_r$. To conclude we mainly need to prove that $\psi(t)$ does converge to $\psi_\infty$ with rate $\rho(\theta_\rho\,t)$ for some positive $\theta_\rho$. In order to do so, we write the Rankine-Hugoniot condition as 
\begin{align*}
\psi'(t)\,&=F(\uu_r(\psi(t)),\uu_\ell(\psi(t)-\psi_{char})\,)
-F(\uu_r(\psi_\infty),\uu_\ell(\psi_\infty-\psi_{char})\,)\\
&\quad+\,F(u_r(t,\sigma\,t+\psi(t)),u_\ell(t,\sigma\,t+\psi(t))\,)
-F(\uu_r(\psi(t)),\uu_\ell(\psi(t)-\psi_{char})\,)\,,
\end{align*}
and deduce (under a smallness condition on $E_0$ independent of $\rho$)
\[
|\psi(t)-\psi_\infty|\leq 
K_{RH}''E_0\,\eD^{-\omega\,t}\,
+K_\rho'E_0^{[\rho]}\,\int_0^t \eD^{-\omega\,(t-s)} \rho(\theta_\rho''\,s)\,\dD s\,,
\]
for some positive $\omega$, $K_{RH}''$, $K_\rho'$, $\theta_\rho''$. This is sufficient to derive the claimed bound on $\psi(t)-\psi_\infty$ and conclude the proof.
\end{proof}

\section{Multidimensional equations}
\label{s:multiD} 

\subsection{Framework and classification}

We now move to the study of multidimensional equations of the form
\begin{equation}
\label{e:multiD}
\partial_t u + \Div_{\bfx} (\bF(u)) = g(u)
\end{equation}
where $u(t,\bfx)\in\R$, $\bfx\in \R^d$, $d\in\N$ and $\bF:\R\to\R^d$.

We are interested in the stability of planar traveling waves, that is, in solutions of the form $(t,\bfx)\mapsto \uu(\bfeta\cdot(\bfx-t\,\bfsigma))$ for some direction $\bfeta\in\bS^{d-1}$ and some wave speed $\bfsigma\in\R^d$. Note that at this stage only the component of $\bfsigma$ in the direction $\bfeta$ plays a role.

Since we shall not need to modulate direction parameters\footnote{We refer to \cite{Audiard-Rodrigues} for an example where in contrast this plays a deep role.}, without loss of generality we may assume that $\bfeta$ points in the direction of the first coordinate. Accordingly we split $\bfx=(x,\bfy)$, $\bF=(\fpar,\fperp)$, $\bfsigma=(\sigpar,\sigperp)$ and we restrict to waves of the form $(t,\bfx)\mapsto \uu(x-t\,\sigpar)$. Then Assumption~\ref{as:generic} is understood as holding with $\fpar$ replacing $f$.

When dealing with such a kind of waves, space-modulated and orbital stabilities require controls respectively of $\delta_{sm}(u(t,\cdot+\bfsigma\,t),\uu)$ and $\delta_{orb}(u(t,\cdot+\bfsigma\,t),\uu)$ where\footnote{We denote as $\uu$ both $x\mapsto \uu(x)$ and $\bfx\mapsto \uu(x)$.}
\begin{align*}
\delta_{sm}(v_1,v_2)&:=\hspace{-2.25em}
\inf_{\substack{\Psi:\,\R^d\to\R\textrm{ such that for any }\bfy\\\Psi(\cdot,\bfy)\ \textrm{is a}\ \cC^1\ \textrm{diffeomorphsim}}}
\hspace{-2.25em}
(\|\bfx\mapsto v_1(\Psi(x,\bfy),\bfy)-v_2(\bfx)\|_{W^{1,\infty}_\Omega((\R\setminus D)\times\R^{d-1})}
+\|\d_x\Psi-1\|_{W^{1,\infty}_\Omega(\R^d)})\,\\[1em]
\delta_{orb}(v_1,v_2)&:=\inf_{\varphi:\,\R^{d-1}\to\R}
\|\bfx \mapsto v_1(x+\varphi(\bfy),\bfy)-v_2(\bfx)\|_{W^{1,\infty}_\Omega((\R\setminus D)\times\R^{d-1})}\,.
\end{align*}
Moreover the weight $\Omega$ is constrained to be constant equal to one for classical weightless stability and to be a function of $x$ only for convective stability.

\bt\label{th:class-multiD}
Under Assumption~\ref{as:generic}, for any one-dimensional non-degenerate wave of \eqref{e:original-eq} in the sense of Definition~\ref{def:non-degenerate} one-dimensional convective stability implies multidimensional convective nonlinear stability and one-dimensional classical weightless stability implies multidimensional classical weightless nonlinear stability. Moreover for stable waves with no characteristic point and no discontinuity between two non-constant portions nonlinear asymptotic stability holds (either convective or classical weightless) whereas for other stable waves any solution arising from a small perturbation converges to a possibly genuinely multidimensional nearby wave.
\et

The part of the latter theorem devoted to constant states or Riemann shocks that are weightlessly stable is already contained in \cite{DR1}. However the multidimensional part of \cite{DR2} devoted to other weightlessly stable waves is restricted to perturbations supported away from characteristic points. This severe restriction drastically simplifies the analysis and the dynamics since it implies that no multidimensional wave is needed to capture the large-time asymptotic behavior. We regard the elucidation of the latter as one of our main contributions.

Let us point out that when proving convergence to a multidimensional wave we use a genuinely multidimensional version of space-modulated stability, controlling 
\begin{align*}
\delta_{sm}(v_1,v_2)&:=\hspace{-1em}
\inf_{\Psi\ \cC^1\ \textrm{diffeomorphsim}}
\hspace{-1em}
(\|v_1\circ \Psi-v_2\|_{W^{1,\infty}_\Omega(\R^d\setminus \cD)}
+\|\nabla\Psi-\I\|_{W^{1,\infty}_\Omega(\R^d)})\,
\end{align*} 
with $v_1=u(t,\cdot+\bfsigma\,t)$ and $v_2$ a wave profile of discontinuity set $\cD$ and speed $\bfsigma$.

\subsection{Nearby multidimensional waves}

To begin with, we discuss the structure of multidimensional waves that are near convectively stable planar waves in a space-modulated sense. We only investigate the existence part here but the asymptotic behavior results do bring the associated uniqueness. 

\subsubsection*{Smooth waves}

To begin with, we consider a smooth non constant convectively stable planar wave of profile $\uu$ and longitudinal speed $\sigpar$. The choice of the transverse part of the speed $\sigperp$ is precisely determined to enforce the existence of nearby multidimensional waves with the same speed $\bfsigma$. Consistently with the definition of $\delta_{sm}$, we search for a wave $(t,\bfx)\mapsto \cU(\bfx-t\,\bfsigma)$ with profile
\[
\cU(\bfx)\,=\,\uu(\Phi(\bfx))
\]
where 
\begin{align*}
\Phi(\bfx)&=(\Id_\R+\psi(\cdot,\bfy))^{-1}(x)\,,&
\d_x\psi&\textrm{ small}\,.
\end{align*}
The smallness of $\d_x\psi$ ensures the invertibility property needed to define $\Phi$. Plugging this \emph{ansatz} in \eqref{e:multiD} one derives that it does provide a traveling wave solution if and only if
\[
(\fpar'(\uu)-\sigpar)\,\d_x\psi+(\fperp'(\uu)-\sigperp)\cdot\nabla_\bfy\psi\,=\,0\,.
\]
This may be solved with 
\begin{align*}
\psi(\bfx)&=\psi_0(\bfy-\bZ(x))\,,&
(\fpar'(\uu)-\sigpar)\,\bZ'&=\fperp'(\uu)-\sigperp\,,
\end{align*}
where $\psi_0$ is arbitrary. When $\uu$ takes a characteristic value $u_\star$ we need to choose $\sigperp=\fperp'(u_\star)$ so as to solve the $\bZ$-equation. For such waves, in view of the one-dimensional stability results it is natural to choose $\bZ(x_\star)=0$, hence
\begin{align*}
\bZ(x)&\,=\,\int_{x_\star}^x\frac{\fperp'(\uu)-\sigperp}{\fpar'(\uu)-\sigpar}
=\int_{u_\star}^{\uu(x)}\frac{\fperp'-\sigperp}{g}\,,&
\sigperp&\,=\,\fperp'(u_\star)\,,
\end{align*}
where $x_\star$ is the characteristic point of $\uu$, so that $x=x_\star+\psi_0(\bfy)$ is the equation of the characteristic hypersurface of the multidimensional pattern. 

There remains to consider smooth monostable fronts. For the sake of definiteness, we assume, without loss of generality, that the unstable endstate is at $+\infty$ and denote it as $u_{+\infty}$. For those, in view of the one-dimensional stability results it is natural to choose $\bZ(+\infty)=0$, hence
\begin{align*}
\bZ(x)&\,=\,-\int_x^\infty\frac{\fperp'(\uu)-\sigperp}{\fpar'(\uu)-\sigpar}
=-\int_{\uu(x)}^{u_{+\infty}}\frac{\fperp'-\sigperp}{g}\,,&
\sigperp&\,=\,\fperp'(u_{+\infty})\,,
\end{align*}
so that 
\[
\lim_{x\to\infty}\eD^{\frac{g'(u_{+\infty})}{|\fpar'(u_{+\infty})-\sigpar|}\,x}(\cU(x,\bfy)-\uu(x-\psi_0(\bfy)))\,=\,0\,.
\]

Note that, whatever the kind of $\uu$, the family of multidimensional waves is parameterized by an infinite-dimensional object, the shift function $\psi_0$, and each member of the family is a genuinely\footnote{Except when both $\fperp$ and $\psi_0$ are affine functions.} multidimensional wave, not a one-dimensional wave in disguise.

\subsubsection*{Discontinuous waves}

We now move to discontinuous waves. 

We begin with Riemann shocks, whose treatment is straightforward. Let $u_{-\infty}$ and $u_{+\infty}$ denote the endstates of a Riemann shock $\uu$ with longitudinal speed $\sigpar$ and discontinuity at $0$. For any $\cC^1$ $\psi_0:\R^{d-1}\to\R$ with small gradient, one obtains a traveling wave solution to \eqref{e:multiD} of profile $\cU$ and speed $\bfsigma$ by setting
\begin{align*}
\sigperp&\,=\,\frac{[\fperp(\uu)]_0}{[\uu]_0}
\,=\,\frac{\fperp(u_{+\infty})-\fperp(u_{-\infty})}{u_{+\infty}-u_{-\infty}}\,,\\
\cU(\bfx)&\,=\,
\begin{cases}
u_{-\infty}&\quad\textrm{if }x<\psi_0(\bfy)\,,\\
u_{+\infty}&\quad\textrm{if }x>\psi_0(\bfy)\,.
\end{cases}
\end{align*}
Note that the smallness of $\nabla_\bfy \psi_0$ is used only to check that Oleinik's conditions persist. For a sharper description of the existence region for this kind of piecewise-constant wave we refer the reader to \cite{Serre_scalar} whose existence part is readily adapted to the present case.

Now we explain how to glue together two parts of profiles of the same speed $\bfsigma$, at least one being non constant. The discussion is analogous to the one preceding Theorem~\ref{th:mixed}. 

Let $\uu$ and $\sigpar$ define a one-dimensional non-degenerate traveling wave with discontinuity at $0$ such that $[g(\uu)]_0/[\uu]_0<0$. Then define one-dimensional extensions $\uu_\ell$ and $\uu_r$ as in the discussion preceding Theorem~\ref{th:mixed}. Now let $\sigperp$, $\varphi_\ell$ and $\varphi_r$ be such that
\begin{align*}
\cU_\ell(\bfx)&:=\uu_\ell(x+\varphi_\ell(\bfx))\,,&
\cU_r(\bfx)&:=\uu_r(x+\varphi_r(\bfx))\,,&
\end{align*} 
define $\cU_\ell$ and $\cU_r$ satisfying equations for multidimensional profiles of speed $\bfsigma$ on neighborhoods respectively of $(-\infty,0]\times \R^{d-1}$ and $[0,\infty)\times \R^{d-1}$, with $(\varphi_\ell,\varphi_r)$ small with small $x$-derivative. We want to show that there exists $\psi_0$ small such that
\[
\R^d\setminus\{\,(\psi_0(\bfy),\bfy)\,;\,\bfy\in\R^{d-1}\,\}\to\R\,,\qquad x\mapsto\begin{cases}
\cU_\ell(\bfx)&\quad \textrm{if }x<\psi_0(\bfy)\,,\\
\cU_r(\bfx)&\quad \textrm{if }x>\psi_0(\bfy)\,,
\end{cases}
\]
defines a wave profile of speed $\bfsigma$. For such an \emph{ansatz}, the latter requirement reduces to the Rankine-Hugoniot conditions and is equivalent to 
\begin{align*}
\Big(\fperp(\cU_r)-\sigperp\cU_r&-\left(\fperp(\cU_\ell)-\sigperp\cU_\ell\right)\Big)(\psi_0(\bfy),\bfy)\cdot \nabla_\bfy\psi_0(\bfy)\\
&\,=\,\Big(\fpar(\cU_r)-\sigpar\cU_r-\left(\fpar(\cU_\ell)-\sigpar\cU_\ell\right)\Big)(\psi_0(\bfy),\bfy)\,,\qquad\bfy\in\R^{d-1}\,.
\end{align*}

The proof that a solution $\psi_0$ does exist, that it is unique (among small solutions) and that it depends continuously on $(\varphi_\ell,\varphi_r)$ does not follow readily from an implicit function theorem. The existence part may be obtained by proving that solutions to
\begin{align*}
[\uu]_0\,\d_t\psi(t,\bfy)
+\Big(\fperp(\cU_r)&-\sigperp\cU_r-\left(\fperp(\cU_\ell)-\sigperp\cU_\ell\right)\Big)(\psi(t,\bfy),\bfy)\cdot \nabla_\bfy\psi(t,\bfy)\\
&\,=\,\Big(\fpar(\cU_r)-\sigpar\cU_r-\left(\fpar(\cU_\ell)-\sigpar\cU_\ell\right)\Big)(\psi(t,\bfy),\bfy)\,,\qquad\bfy\in\R^{d-1}\,,
\end{align*}
starting with small initial data possess a limit as $t\to\infty$. The key observation to prove the latter is that the difference between two solutions decays exponentially in time, as follows from $[g(\uu)]_0/[\uu]_0<0$ provided that $(\varphi_\ell,\varphi_r)$ is sufficiently small. The latter observation also yields the claimed uniqueness for $\psi_0$. The sharp continuity property requires a little more care. We omit the details of the arguments since those are both classical and cumbersome.

\subsection{Stability of smooth characteristic fronts}

We provide detailed versions of Theorem~\ref{th:class-multiD} only for smooth characteristic fronts. With this analysis in hands, the rest of the cases is concluded along the gluing strategy already used in the preceding section and we refer to \cite{DR1} for a multidimensional instance of the latter.

The analysis being inherently anisotropic we introduce relevant anisotropic norms $\|\cdot\|_{W^{1,\infty}_\shortparallel}$, $\|\cdot\|_{W^{2,\infty}_\shortparallel}$ defined by
\begin{align*}
\|v\|_{W^{1,\infty}_\shortparallel}
&:=\|\,(v,\,\d_x v)\|_{L^\infty(\R^d)}\,,&
\|v\|_{W^{2,\infty}_\shortparallel}
&:=\|\,(v,\,\nabla_\bfx v)\|_{W^{1,\infty}_\shortparallel}\,.
\end{align*}
With any weight $\rho:\R_+\to\R_+$ we associate weights $\rho_\shortparallel:\,\R^d\to\R_+$ and  $\rho_{-1}:\,\R^*\times\R^{d-1}\to\R_+$ through
\begin{align*}
\rho_\shortparallel(\bfx)
&:=\rho(|x|)\,,&
\rho_{-1}(\bfx)
&:=\frac{\rho(|x|)}{|x|}\,,
\end{align*}
and introduce corresponding norms $\|\cdot\|_{L^{\infty}_{\rho,\shortparallel}}$, $\|\cdot\|_{L^{\infty}_{\rho,-1}}$ and $\|\cdot\|_{W^{1,\infty}_{\rho,\shortparallel}}$ by
\begin{align*}
\|v\|_{L^{\infty}_{\rho,\shortparallel}}
&:=\|\,\rho_\shortparallel\,v\|_{L^\infty(\R^d)}\,,&
\|v\|_{L^{\infty}_{\rho,-1}}
&:=\|\,\rho_{-1}\,v\|_{L^\infty(\R^d)}\,,&
\|v\|_{W^{1,\infty}_{\rho,\shortparallel}}
&:=\|\,\rho_\shortparallel\,(v,\,\d_x v)\|_{L^\infty(\R^d)}\,.
\end{align*}
It is important to understand that here the weights are used to let the functions grow, whereas in the convective stability results they are used to enforce decay. In particular the extra linear growth encoded by $\|\cdot\|_{L^{\infty}_{\rho,-1}}$ originates in the fact that the $\bZ$-functions of the previous subsection typically grow linearly.

\bt\label{th:multiD}
Consider a one-dimensional non-degenerate traveling wave of profile $\uu$ and speed $\sigpar$ such that $\uu$ is smooth with a single characteristic point at $0$ and  with $g'(\uu(0))>0$. In particular, denoting as $u_-$ and $u_+$ the corresponding endstates, $g'(u_-)<0$ and $g'(u_+)<0$.\\
Set $u_\star:=\uu(0)$, 
\begin{align*}
\sigperp\,&:=\,\fperp'(u_\star)\,,&
&\textrm{and}&
\omega_0
&:=\min(\{g'(u_\star),|g'(u_-)|,|g'(u_+)|\})\,.
\end{align*}

There exist positive $\delta_0$ and $C$ and for any sub-exponential weight $\rho$ there exist positive $C_\rho$ and $\theta_\rho$ such that if $u_0\in BUC^2(\R^d)$ and
\[
E_0:=\|u_0-\uu\|_{W^{2,\infty}(\R^d)}
\leq \delta_0
\]
then there exists a $\cC^1$ phase shift $\psi:\,\R_+\times\R^d\to \R$ such that the solution $u$ to \eqref{e:multiD} starting from $u(0,\cdot)=u_0$ satisfies for any $t\geq0$ 
\begin{align*}
\|u(t,\cdot+\bfsigma\,t+(\psi(t,\cdot),0))-\uu\|_{W^{2,\infty}_\shortparallel}
&\leq\,C\,\eD^{-\omega_0\,t}\,E_0\,\\
\|(\psi(t,\cdot),\d_t\psi(t,\cdot),\nabla_\bfx\psi(t,\cdot))\|_{L^\infty}
&\leq\,C\,E_0\,
\end{align*}
and, denoting $\cU$ the unique multidimensional profile near $\uu$ such that $\cU^{-1}(\{\,u_\star\,\})\,=\,u_0^{-1}(\{\,u_\star\,\})$, there exists a $\cC^1$ phase shift $\varphi:\,\R_+\times\R^d\to \R^{d-1}$ such that for any $t\geq0$ 
\begin{align*}
\|u(t,\cdot+\bfsigma\,t+(0,\varphi(t,\cdot)))-\cU\|_{W^{2,\infty}_\shortparallel}
&\leq\,C\,\eD^{-\omega_0\,t}\,E_0\,
\end{align*}
and
\begin{enumerate}
\item if $\Omega\in L_{loc}^\infty(\R_+;\R_+)$ is such that $\lim_{+\infty} \Omega=0$ then
\[
\lim_{t\to\infty}
(\|\,\varphi(t,\cdot)\|_{L^{\infty}_{\Omega,-1}}
+\|\,(\d_t\varphi(t,\cdot),\nabla_\bfx \varphi(t,\cdot))\|_{L^{\infty}_{\Omega,\shortparallel}})\,=\,0\,;
\]
\item if $\rho$ is a sub-exponential weight then 
\begin{align*}
\|\varphi(t,\cdot)\|_{L^\infty_{\rho,-1}}
+\|(\d_t\varphi(t,\cdot),\nabla_\bfx \varphi(t,\cdot))\|_{W^{1,\infty}_{\rho,\shortparallel}}
\,\leq\,C_\rho\,E_0\,\rho(\,\theta_\rho\,t)\,.
\end{align*}
\end{enumerate}
\et

\begin{proof}
The first key point is to notice that $u$ may be recovered through
\[
u(t,x+t\,\sigpar,\bY(t,x,\bfy)+t\,\sigperp)
\,=\,v(t,x,\bfy)
\]
from $v$ and $\bY$ solving
\begin{align}\label{eq:shape}
\d_tv+(\fpar'(v)-\sigpar)\d_xv&=\,g(v)\\
\label{eq:position}
\d_t\bY+(\fpar'(v)-\sigpar)\d_x\bY&=\,\fperp'(v)-\sigperp
\end{align}
with initial data
\begin{align*}
v(0,\cdot)&=\,u_0\,,&
\bY(0,x,\bfy)&=\,\bfy\,.
\end{align*}
Equation~\eqref{eq:shape} is almost directly the one-dimensional equation studied in \cite[Section~4]{DR2} whereas System~\eqref{eq:position} is affine in $\bY$ and essentially one-dimensional.

To untangle the large-time dynamics of $v$, we introduce $\psi_{char}:\R^{d-1}\to\R$ such that 
\[
u_0(\psi_{char}(\bfy),\bfy)\,=\,u_\star\,,
\]
and set 
\[
w(t,x,\bfy):=v(t,x+\psi_{char}(\bfy),\bfy)\,.
\]
Note that $w(0,0,\cdot)\equiv u_\star$ and
\[
\d_tw+(\fpar'(w)-\sigpar)\d_xw\,=\,g(w)\,.
\]
Provided that $E_0$ is sufficiently small, we deduce from \cite[Section~4]{DR2} that
\begin{align*}
\|w(t,\cdot)-\uu\|_{W^{1,\infty}_\shortparallel}
&\leq K\,\eD^{-\omega_0\,t}\,\|w(0,\cdot)-\uu\|_{W^{1,\infty}_\shortparallel}
\leq K'\,\eD^{-\omega_0\,t}\,E_0\\
\|\nabla_\bfx(w(t,\cdot)-\uu)\|_{W^{1,\infty}_\shortparallel}
&\leq K\,\eD^{-\omega_0\,t}\,\|\nabla_\bfx(w(0,\cdot)-\uu)\|_{W^{1,\infty}_\shortparallel}
\leq K'\,\eD^{-\omega_0\,t}\,E_0\,,
\end{align*}
for some $K$, $K'$.

To unravel the asymptotically one-dimensional structure of Equation~\ref{eq:position} we introduce 
\begin{align*}
\bY_{char}(t,x,\bfy)&:=\bY(t,x+\psi_{char}(\bfy),\bfy)\,,&
\end{align*}
so that 
\begin{align*}
\d_t\bY_{char}+(\fpar'(w)-\sigpar)\d_x\bY_{char}&=\,\fperp'(w)-\sigperp\,,&
\bY_{char}(0,x,\bfy)&=\,\bfy\,.
\end{align*}
For the sake of comparison we also introduce $\bY_\infty(x,\bfy):=\bfy+\bZ(x)$ where $\bZ$ is as in the previous subsection
\begin{align*}
\bY_\infty(x,\bfy)&:=\bfy+\bZ(x)\,,&
\bZ(x)&\,=\,\int_{x_\star}^x\frac{\fperp'(\uu)-\sigperp}{\fpar'(\uu)-\sigpar}
=\int_{u_\star}^{\uu(x)}\frac{\fperp'-\sigperp}{g}\,.
\end{align*}
We recall that then
\[
\cU(x+\psi_{char}(\bfy),\bY_\infty(x,\bfy))=\uu(x)\,.
\]
Accordingly in the end we shall set 
\[
\psi(t,x,\cdot):=\psi_{char}\circ \bY_{char}(t,x,\cdot)^{-1}\,.
\]
and
\begin{align*}
\cY(t,x,\cdot)&:=
\bY_{char}(t,x-\psi_{char}(\cdot),\cdot)
\circ \left(\bY_\infty(x-\psi_{char}(\cdot),\cdot)\right)^{-1}\,,&
\varphi(t,\bfx)&:=\cY(t,\bfx)-\bfy\,.
\end{align*}

Note that
\begin{align*}
\d_t(\bY_{char}-\bY_\infty)+(\fpar'(w)-\sigpar)\d_x(\bY_{char}-\bY_\infty)
&=\,\fperp'(w)-\fperp'(\uu)-(\fpar'(w)-\fpar'(\uu))\d_x\bY_\infty\,,&\\
(\bY_{char}-\bY_\infty)(0,x,\bfy)&=\,-\bZ(x)\,.
\end{align*}
Differentiating the latter and combining with bounds on $w-\uu$ yield for some $K''$
\begin{align*}
\|\d_x(\bY_{char}-\bY_\infty)(t,\cdot)\|_{L^\infty}
&\leq\,K''\,(1+E_0)\,,&\\
\|\nabla_\bfy(\bY_{char}-\bY_\infty)(t,\cdot)\|_{L^\infty}
&\leq\,K''\,E_0\,.
\end{align*}
This is sufficient to justify the invertibility of $\bY_{char}$ when $E_0$ is small and derive the claimed comparison with $\uu$ and bounds on $\psi$.

The end of the proof uses more crucially that
\[
(\bY_{char}-\bY_\infty)(0,0,\cdot)\,\equiv\, 0\,,
\]
and follows readily from the next proposition and the observation already used and proved along the proof of Theorem~\ref{th:fronts} that if $\lim_{+\infty}\Omega=0$ then $\Omega$ is bounded from above by a sub-exponential weight.
\end{proof}

We recall from the proof of Proposition~\ref{p:rho} that $\cT_a(t_1,t_2)$ denotes the solution operator from time $t_1$ to $t_2$ associated with
\[
\d_t w+a\,\d_x w\,=\,0
\]
We see the following proposition as somehow dual to Proposition~\ref{p:rho}.

\bpr\label{p:unrho}
Under the assumptions of Theorem~\ref{th:multiD}, 
there exists a positive $\delta$ and for any sub-exponential weight $\rho$, there exist positive $\eta_\rho$ and $K_\rho$ such that if $T>0$ is such that $a(\cdot,0,\cdot)\equiv 0$ and
\[
\|a-(f'(\uu)-\sigma)\|_{L^\infty(0,T;W^{1,\infty}_\shortparallel)}\leq \delta
\]
then for any $0\leq s\leq t\leq T$, any sub-exponential weight $\rho$ and any $0<\eta\leq \eta_\rho$
\begin{align*}
\|\d_x\cT_a(s,t)w\|_{L^{\infty}_{\rho,\shortparallel}}
&\leq K_\rho\,\rho(\eta\,(t-s))\,\|\d_x w\|_{L^{\infty}}
\exp\left(K_\rho\,\int_s^t\|a(\tau,\cdot)-(f'(\uu)-\sigma)\|_{W^{1,\infty}_{\shortparallel}}\dD\tau\right)\,,
\end{align*}
and if $w(0,\cdot)\equiv 0$
\begin{align*}
\|\cT_a(s,t)w\|_{W^{1,\infty}_{\rho,\shortparallel}}
&\leq K_\rho\,\rho(\eta\,(t-s))\,\|w\|_{W^{1,\infty}_{\shortparallel}}
\exp\left(K_\rho\,\int_s^t\|a(\tau,\cdot)-(f'(\uu)-\sigma)\|_{W^{1,\infty}_{\shortparallel}}\dD\tau\right)\,,\\
\|\cT_a(s,t)w\|_{L^{\infty}_{\rho,-1}}
&\leq K_\rho\,\rho(\eta\,(t-s))\,\|\d_xw\|_{L^{\infty}}
\exp\left(K_\rho\,\int_s^t\|a(\tau,\cdot)-(f'(\uu)-\sigma)\|_{W^{1,\infty}_{\shortparallel}}\dD\tau\right)\,.
\end{align*}
\epr

\begin{proof}
We provide details only for the $W^{1,\infty}_{\rho,\shortparallel}$ bound. The $L^{\infty}_{\rho,\shortparallel}$ on the derivative is proved similarly and the $L^{\infty}_{\rho,-1}$ bound is then deduced from it by integration.

The arguments of \cite[Section~4]{DR2} are readily adapted to yield 
\begin{align*}
\|\cT_a(s,t)w\|_{W^{1,\infty}_{\shortparallel}}
&\leq C\,\|w\|_{W^{1,\infty}_{\shortparallel}}
\exp\left(C\,\int_s^t\|a(\tau,\cdot)-(f'(\uu)-\sigma)\|_{W^{1,\infty}_{\shortparallel}}\dD\tau\right)\,,
\end{align*}
and
\begin{align*}
&\|\cT_a(s,t)w\|_{W^{1,\infty}_{\shortparallel}([-M,M]\times \R^{d-1})}\\
&\quad\leq C_M\,\eD^{-g'(u_\star)\,(t-s)}\,\|w\|_{W^{1,\infty}_{\shortparallel}([-M,M]\times \R^{d-1})}
\exp\left(C_M\,\int_s^t\|a(\tau,\cdot)-(f'(\uu)-\sigma)\|_{W^{1,\infty}_{\shortparallel}}\dD\tau\right)\,,
\end{align*}
for some constants $C$, $C_M$, with $M>0$ being arbitrary. These are sufficiently to deduce the claimed bound on $(\cT_a(s,t)w,\d_x \cT_a(s,t)w)(\bfx)$ respectively when $|x|\geq \theta_0 (t-s)$ for some $\theta_0>0$ and when $|x|\leq M$.

To complete the argument we need to take into account support propagation. For definiteness choose $M=1$ and denote $\theta_\cT>0$ a lower bound for $|a|$ on $(0,T)\times \R\setminus[-1,1]\times \R^{d-1}$. Then if $\theta_0<\theta_\cT$, one obtains likewise
\begin{align*}
\|\cT_a(s,t)w\|_{W^{1,\infty}_{\shortparallel}([-\theta_0\,(t-s),\theta_0\,(t-s)]\times \R^{d-1})}
&\leq C\,\|\cT_a(s,s+(1-\tfrac{\theta_0}{\theta_\cT})(t-s))w\|_{W^{1,\infty}_{\shortparallel}([-M,M]\times \R^{d-1})}\\
&\times\exp\left(C\,\int_{s+(1-\tfrac{\theta_0}{\theta_\cT})(t-s)}^t\|a(\tau,\cdot)-(f'(\uu)-\sigma)\|_{W^{1,\infty}_{\shortparallel}}\dD\tau\right)
\end{align*}
which by combination with the above gives 
\begin{align*}
&\|\cT_a(s,t)w\|_{W^{1,\infty}_{\shortparallel}([-\theta_0\,(t-s),\theta_0\,(t-s)]\times \R^{d-1})}\\
&\quad\leq C'\,\eD^{-g'(u_\star)\,\left(1-\tfrac{\theta_0}{\theta_\cT}\right)\,(t-s)}\,\|w\|_{W^{1,\infty}_{\shortparallel}}\,\exp\left(C'\,\int_s^t\|a(\tau,\cdot)-(f'(\uu)-\sigma)\|_{W^{1,\infty}_{\shortparallel}}\dD\tau\right)
\end{align*}
for some $C'$. Hence the result.
\end{proof}

\bibliographystyle{alphaabbr}
\bibliography{Balance-laws}

\end{document}